\documentclass[a4paper,11pt]{article}
\usepackage[utf8]{inputenc}

\usepackage[titletoc]{appendix} 
\usepackage{xcolor}
\usepackage[pdfencoding=auto, psdextra]{hyperref}
\hypersetup{
    colorlinks=true,
    citecolor=blue
    }
\usepackage{bookmark}
\usepackage{enumitem} 
\usepackage{authblk} 
\usepackage{pifont} 

\usepackage{amssymb}
\usepackage{amsmath}
\usepackage{amsthm}
\usepackage{bm} 
\usepackage{mathtools}
\usepackage{stmaryrd} 
\usepackage{setspace} 
\usepackage{mathabx}
\usepackage{cleveref} 

\crefrangeformat{subequation}%
      {(#3#1#4--#5\crefstripprefix{#1}{#2}#6)}

\theoremstyle{plain}
\newtheorem{theorem}{Theorem}[section]
\newtheorem{proposition}[theorem]{Proposition}
\newtheorem{lemma}[theorem]{Lemma}

\newtheorem{assumption}[theorem]{Assumption}

\newtheorem{remark}[theorem]{Remark}

\usepackage{stmaryrd}
\usepackage{graphicx}

\newcommand{\N}{\mathbb{N}}
\newcommand{\R}{\mathbb{R}}

\def\P{\mathbb{P}}
\newcommand{\E}{\mathbb{E}}

\newcommand{\D}{\mathbb{D}}

\newcommand{\e}{\mathrm{e}}

\newcommand{\Lone}[1]{\left\lVert#1\right\rVert_{1}}
\newcommand{\C}{\mathcal{C}} 
\newcommand{\ind}{\mathbf{1}}
\newcommand{\setind}[1]{\mathbf{1}_{\left\{#1\right\}}}
\newcommand{\bbrackets}[2]{\llbracket #1, #2 \rrbracket}
\newcommand{\angles}[2]{\langle #1, #2 \rangle}

\newcommand{\X}{\mathcal{X}}
\newcommand{\Z}{\mathcal{Z}}
\newcommand{\U}{\mathcal{U}}
\newcommand{\Q}{\mathcal{Q}}
\newcommand{\A}{\mathcal{A}}

\newcommand{\bfz}{\mathbf{z}}
\newcommand{\bfk}{\mathbf{k}}
\newcommand{\bfe}{\mathbf{e}}
\newcommand{\bfx}{\mathbf{x}}
\newcommand{\bfv}{\mathbf{v}}

\newcommand{\configS}{\mathcal{S}}
\newcommand{\generation}[1]{\mathbb{G}(#1)}
\newcommand{\popE}{\E_{\bfz}}
\newcommand{\spineE}{\E_{x, \bfz}}

\newcommand{\couple}{\mathbb{X}}
\newcommand{\couplespine}{\mathbb{Y}^{(t)}}
\newcommand{\Ybb}{\mathbb{Y}}
\newcommand{\calR}{\mathcal{R}}
\newcommand{\calP}{\mathcal{P}}
\newcommand{\proj}[1]{\left.#1\right|_{\X}}
\newcommand{\bbY}{\mathbb{Y}}
\newcommand{\m}{\mathfrak{m}}

\title{Empirical distribution of ancestral lineages in populations with density-dependent interactions}
\author[1,2]{Madeleine Kubasch}
\affil[1]{\footnotesize Centre de Mathématiques Appliquées (CMAP), Ecole Polytechnique, Palaiseau, France}
\affil[2]{\footnotesize Institute of Ecology and Environmental Sciences (iEES), Sorbonne Université, Paris, France}
\date{\today}

\begin{document}

\title{Empirical distribution of ancestral lineages in populations with density-dependent interactions} 
\maketitle

\begin{abstract}
We study a density-dependent Markov jump process describing a population where each individual is characterized by a type, and reproduces at rates depending both on its type and on the population type distribution. 
First, we exhibit a time-inhomogeneous Markov process, which allows to capture the behavior of a sampled lineage in the population process through a many-to-one formula. 
Second, under classical assumptions, the population type distribution converges to a deterministic limit as population size grows to infinity.
Here, we focus on the empirical distribution of ancestral lineages in this large population limit, for which we establish an asymptotic many-to-one formula and quantify the speed of convergence.
\smallskip

\noindent    \textbf{Keywords.} Interactions; Markov jump process; population process; many-to-one formula; large population limit.
\end{abstract}

\section{Introduction}

When considering population processes arising in various fields such as population genetics or epidemiology, the study of ancestral lineages may provide crucial information. For instance, such lineages yield insight on epidemic spread through contamination chains \cite{duchampsGeneralEpidemiologicalModels2023}, or on the evolution of a trait of interest under selection \cite{calvezDynamicsLineagesAdaptation2022}. As a consequence, several methods have been developed to finely characterize those lineages. On the one hand, a classical approach is to consider a backward-in-time process which reconstructs the genealogy by moving back from time $t$ to time $0$, and which is related to the initial population process by duality \cite{berestyckiRecentProgressCoalescent2009, baakeProbabilisticViewDeterministic2018, machRecursiveTreeProcesses2020, corderoGeneralSelectionModels2022, duchampsGeneralEpidemiologicalModels2023}. On the other hand, there also exists a forward-in-time approach which relies on a second population process, with one distinguished individual (the \emph{spine}) whose lineage behaves as the lineage of a sampled individual in the original process. 

More precisely, these \emph{spinal constructions} have originally been introduced for branching processes, using an appropriate change in probability. The key to the construction of the spinal process is that the reproduction rates of the spine are biased towards leaving more numerous descendants than other individuals. This leads to the emergence of size-biased distributions, which accurately depict the survivorship bias induced by sampling. Generally speaking, the obtained spinal construction has several strengths. Notably, it allows to establish \emph{many-to-one} formulas (\emph{e.g.} \cite{lyonsConceptualProofsLLogL1995, georgiiSupercriticalMultitypeBranching2003, harrisBranchingBrownianMotion2016, harrisManytofewLemmaMultiple2017}), which are closely related to Feynman-Kac path equations \cite[Sections 1.3 and 1.4.4]{delmoralFeynmanKacFormulae2004}. Many-to-one formulas translate the expected value of a functional evaluated over the lineages in the branching process, into the expectation of the functional evaluated along the spine, whose trajectories are exponentially weighted to capture the growth of the population. If the exponential weight is deterministic, this immediately implies a numerical advantage for computing such averages through Monte-Carlo simulations. Indeed, simulations of the spine are numerically affordable, whereas simulations of the whole genealogical tree in the original branching process can be numerically challenging due to exponential growth \cite{nagelRealisticAgentbasedSimulation2021}.
Also, spinal constructions have proven an effective way of establishing classical key results on branching processes, such as the Kesten Stigum theorem \cite{lyonsConceptualProofsLLogL1995, georgiiSupercriticalMultitypeBranching2003}. More recently, the semi-group associated to the spinal construction has proven a successful tool in the analysis of non-conservative semi-groups, extending its applications beyond branching processes \cite{bansayeErgodicBehaviorNonconservative2020, bansayeNonconservativeHarrisErgodic2022}. 

While many models for population dynamics arising, for instance, in biology and epidemiology do not satisfy \emph{per se} the branching approximation, a classical approach is to consider regimes in which the population process can be well approached by a branching process, using coupling arguments. For example, in epidemiology, it is well-known that at the beginning of an epidemic, the tree of infections can be captured by a branching process which neglects the depletion in susceptible individuals \cite{ballStrongApproximationsEpidemic1995}. Similarly, in order to analyze the lineage of a uniformly sampled individual in a population which is subject to evolution under a changing environment, one may consider the stationary regime \cite{calvezDynamicsLineagesAdaptation2022}. However, such branching approximations are restricted to specific parts of the dynamics of interest only; see for instance \cite{barbourApproximatingEpidemicCurve2013} and \cite{bansayeSharpApproximationHitting2024} for details in the case of epidemic models and invasion processes. 

In order to address this limitation, there have been developments towards capturing the ancestral lineage of a sampled individual, as well as the whole genealogical tree, in populations with interactions. Recently, a spinal construction has been developed for this setting, focusing on multi-type processes with discrete type space \cite{bansayeSpineInteractingPopulations2024}. The general idea consists in biasing the reproduction rates of the process, both along the spine and outside of it, according to a positive function $\psi$ of the reproducing particle's trait $x$ and the population's type composition $\bfz$. Intuitively, $\psi(x, \bfz)$ can be regarded as the individual's reproductive value or long-term fertility. Hence, when the spine reproduces, offspring with higher values of $\psi$ given the population state are favored, while the offspring of individuals outside of the spine are biased towards rendering the population more favorable for the spine. This spinal construction has since been extended to include more general type spaces \cite{medousSpinalConstructionsContinuous2024}. It further has served to study the convergence of genealogies of density-dependent branching processes to the Kingman coalescent, establishing a connection with backward-in-time approaches \cite{andreMomentsDensitydependentBranching2025}.  

In this paper, we focus on the empirical type distribution of ancestral lineages ($\psi = 1$), which arguably corresponds to the most natural and naive sampling strategy. In particular, we aim to capture the underlying survivorship bias. Notably, in the many-to-one formula associated to the aforementioned $\psi$-spine, spinal trajectories are penalized by an exponential weight which is generally stochastic. Thus the interpretation of the survivorship bias embodied by this spinal process is not straightforward, as this penalization needs to be taken into account. In addition, Monte-Carlo estimations of the many-to-one formula become delicate as rare trajectories may have a tremendous impact. As a consequence, we aim at proposing an alternative spinal construction, whose associated many-to-one formula does not require exponential weighting of trajectories. This is achieved by a time-inhomogeneous spinal process, inspired from a similar auxiliary process which captures the lineage of a uniformly sampled individual in a branching process with large initial population  
\cite{marguetUniformSamplingStructured2019}. 

Further, a natural regime to consider is the large population limit. Under classical assumptions, several large population approximations can frequently be established \cite[Chapter 11]{ethierMarkovProcessesCharacterization1986}. Indeed, these approximations may be more prone to mathematical analysis than the original population process. From a numerical point of view, simulating large population approximations is often less expensive than individual-based models for large population sizes, both in terms of computation time and memory. For instance, numerical evaluation of the  $\psi$-spine many-to-one formula requires simulations of the spinal population process, which is computationally expensive for large population sizes. As a consequence, it is natural to develop approximate sampling strategies, which allow to sample directly in those large population approximations. 
Building on the large population limit of the spinal process introduced in \cite{bansayeSpineInteractingPopulations2024}, we establish a many-to-one formula for sampling in the deterministic large population limit with our time-inhomogeneous spinal process. Finally, we quantify the approximation error which is made by sampling in the large population limit instead of the original population process.

This paper is structured as follows. Sections \ref{sec:res-bounded} and \ref{sec:res-lln} present our contributions regarding the survivorship bias in density-dependent population processes, respectively for populations of finite size and in the large population limit. Proofs are provided in Section \ref{sec:proofs}. Finally, Section \ref{sec:discussion} presents a discussion on our results.

\section{Survivorship bias in density-dependent population processes}
\label{sec:res-bounded}

In this section, our aim is to propose a many-to-one formula without stochastic penalization of spinal trajectories, in an effort to gain insight into the survivorship bias. In order to do so, we will introduce a time-inhomogeneous spinal construction. This approach is inspired by \cite{marguetUniformSamplingStructured2019}. While we are particularly interested in the empirical distribution of ancestral lineages, the results obtained in this Section hold for a large class of sampling strategies, and are thus stated in a general setting.  

\subsection{The population process}

\begin{figure}
	\centering
	\includegraphics[width = 0.75\textwidth]{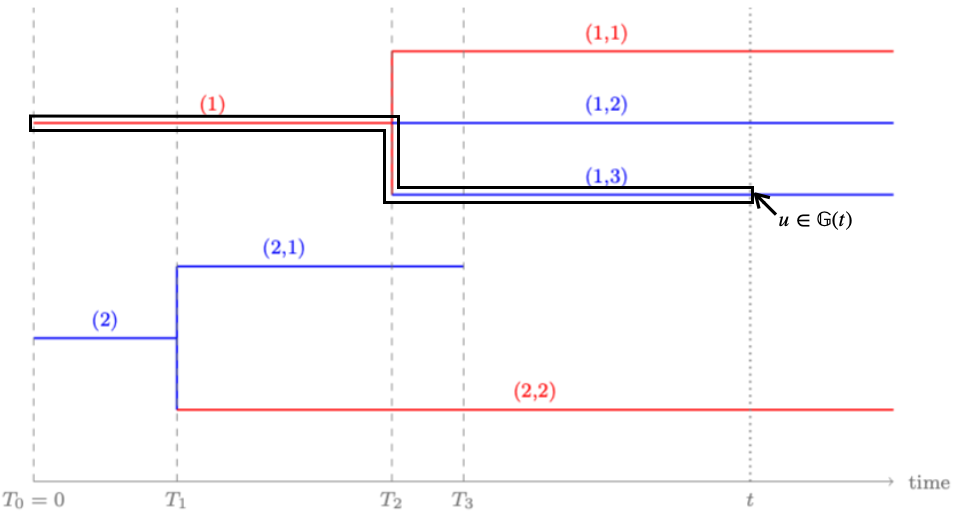}
	\caption{Schematic representation of the population process with type space $\X = \{R, B\}$ designating the colors red ($R$) and blue ($B$). Each individual is represented by a line, and is identified by its label which is indicated above. Reproduction events occur at times $(T_k, k \geq 0)$. At time $t$, the state of the population process is given by $X(t) = \delta_{((1,1), R)} +  \delta_{((1,2), B)} + \delta_{((1,3), B)} + \delta_{((2,2), R)}$. The population contains $Z_R(t) = 2$ individuals of type $R$ and $Z_B(t) = 2$ individuals of type $B$. The ancestral lineage of $u = (1,3) \in \generation{t}$ is framed in black. }
	\label{fig:population-process}
\end{figure}

We consider a structured population, where each individual has a type $x \in \X$, and we assume for convenience that the type space $\X$ is finite. The number of individuals of type $x$ in the population is referred to as $\bfz_x$, and the corresponding vector $\bfz \in \Z = (\N \cup \{0\})^\X$ describes the composition of the population.

Further, individuals will reproduce at rates depending on their type and the current population state. More precisely, an individual of type $x$ may be replaced by some offspring $\bfk = (\bfk_y, y \in \X) \in \Z$, meaning that the individual dies and for any $y \in \X$, $\bfk_y$ individuals of type $y$ are born. This occurs at rate $\tau_\bfk(x,\bfz)$. In addition, we work under the following assumption:

\begin{assumption} 
\label{hyp:avg} 
We assume that the following inequality holds:
\[ \sum_{\bfk \in \Z} \| \bfk \|_1 \| \tau_\bfk \|_{\infty} < +\infty. \]
\end{assumption}

In order to keep track of the genealogy, we make use of the Ulam-Harris-Neveu notations. Let $\U = \bigcup_{k \geq 1} \N^k$, then $u = (u_1, \dots, u_k) \in \U$ represents the $u_k$-th descendant of $(u_1, \dots, u_{k-1})$ and for $u, v \in U$ we write $u \succeq v$ if $v$ is an ancestor of $u$ . The type of $u \in \U$ will be called $x_u$. Hence when an individual $u$ is replaced by its offspring $\bfk$, the new individuals are $(u,1), \dots, (u,\Lone{\bfk})$ and we need to decide the type of each descendant. We thus consider a probability distribution $\Q_\bfk$ on 
\begin{equation*}
	\X_\bfk = \{\bfx \in \X^{\Lone{k}}: \forall x \in \X,  \#\{i: \bfx_i = x\} = \bfk_x\},
\end{equation*}
and $(x_{(u, i)}: i \in \bbrackets{1}{\Lone{k}})$ is distributed as $\Q_\bfk$.

Let us now introduce the stochastic process of interest. Intuitively, it corresponds to describing the set of individuals alive and their types, at each time $t \geq 0$. 
We start from an initial set of individuals $\generation{0} = \mathfrak{g} \subset \N$, and the population will evolve as explained above. At each time $t$, let $\generation{t} \subset \U$ be the set of individuals alive. The process of interest $(X(t), t \geq 0)$ is a Markov jump process with càdlàg trajectories, which for $t \geq 0$ tracks the individuals alive at time $t$, and their types. It can thus be defined informally by 
\[ \forall t \geq 0, \quad X(t) = \sum_{u \in \generation{t}} \delta_{(u, x_u)} \in \mathcal{M}_P(E),\]
where $E = \{(u, x_u) : u \in \U, x_u \in \X \}$ and $\mathcal{M}_P(E)$ is the set of positive point measures on $E$.
In particular, notice that there cannot be explosion, since individuals reproduce at rate less than $\sum_{\bfk \in \Z} \| \tau_\bfk \|_\infty$, which is finite according to Assumption \ref{hyp:avg}. 

Further, $Z(t) = (Z_x(t))_{x \in \X} \in \Z$ yields the composition of the population at time $t$, i.e. for any $x \in \X$,
\[Z_x(t) = \# \{ u \in \generation{t} : x_u = x\}. \] 
Finally, for $u \in \generation{t}$ and $s \leq t$, $x_u(s)$ stands for the type of the unique ancestor of $u$ alive at time $s$. We refer to Figure \ref{fig:population-process} for an illustration of the population process.

Let us end this section by introducing some notation. Define the set 
\begin{equation*}
\configS = \{(x, \bfz) \in \X \times \Z: \bfz_x \geq 1\}.
\end{equation*} 
For any $\bfz \in \Z$, we consider an (arbitrary) labeling $\mathfrak{g}(\bfz) \subset \N$ of individuals, and for $x$ such that $\bfz_x > 0$, we fix $u_x \in \mathfrak{g}(\bfz)$ of type $x$. We denote the corresponding population state by $\mathfrak{X}(\bfz)$. We let $\popE$ and $\P_{\bfz}$ denote the expectation and probability conditionally on $X(0) = \mathfrak{X}(\bfz)$. Considering that an empty sum equals zero,  
\begin{equation*}
	M_t f(x, \bfz) = \popE[\sum_{u \in \generation{t}, u \succeq u_x} f(x_u(t), Z(t))]
\end{equation*}
is a semi-group whose generator $G$ is defined by its action on functions $f: \configS \to \R_+$:
\begin{equation}
\label{eq:defG}
\begin{aligned}
	\forall (x, \bfz) \in \configS_K, \quad & Gf(x, \bfz) = \sum_{\bfk \in \Z} \tau_\bfk(x, \bfz) \Bigg(\sum_{y \in \X} \bfk_y f(y, \bfz + \bfk - \bfe(x)) - f(x, \bfz)\Bigg) \\
	&+ \sum_{\substack{y \in \X \\ \bfk \in \Z}} (\bfz_y - \setind{x = y})\tau_\bfk(y, \bfz) (f(x, \bfz + \bfk - \bfe(y)) - f(x, \bfz)).
\end{aligned}
\end{equation}
We are now ready to turn to our contributions.

\subsection{A many-to-one formula}

Consider a positive function $\psi$ on $\X \times \Z$. As it is kept fixed, dependence on $\psi$ is not explicit in our notations for readability. Define the following function $m$ on $\configS \times [0,t]$:
\begin{equation}
\label{eq:mpsi}
	m(x, \bfz, t) = \popE \Big[\sum_{\substack{u \in \generation{t} \\ u \succeq u_x}} \psi(x_u(t), Z(t)) \Big].
\end{equation}
Notice that, by the Markov property, for $s \in [0,t]$,
\begin{equation*}
	m(x, \bfz, t-s) = \E\left[\sum_{\substack{u \in \generation{t} \\ u \succeq u_x}} \psi(x_u(t), Z(t)) \Bigg| X(s) = {\mathfrak{X}(\bfz)} \right].
\end{equation*}
In words, $m(x, \bfz, t-s)$ corresponds to the $\psi$-weighted average of the types of individuals alive at time $t$ who descend from a given individual of type $x$ at time $s$, given that at time $s$, the population was in state $\bfz$. For instance, if $\psi = 1$, this yields the average number of individuals alive at time $t$, who descend from an individual of type $x$ at time $s$ when the population was in state $\bfz$. 
\medskip

Let us now introduce the time-inhomogeneous spinal process, which allows to capture the behavior of the ancestral lineage of a $\psi$-weighted sample of the population process. For $t \geq 0$ fixed, we will consider the time-inhomogeneous Markov process $(Y^{(t)}(s), \zeta^{(t)}(s))_{s \leq t}$ defined as follows.
The main idea is to follow the type $Y^{(t)}$ of a distinguished individual, which will be referred to as \emph{spine} in analogy to classical spinal constructions. At time $s \leq t$, when of type $x$ in a population of state $\bfz$, the spine is replaced by offspring $\bfk$ and switches to type $y$ (\emph{i.e.} among its offspring $\bfk$, one individual of type $y$ is chosen to become the new spine) with rate 
\begin{equation}
\label{eq:rho}
	\rho^{(t)}_{y, \bfk}(s,x,\bfz) = \tau_\bfk(x, \bfz) \bfk_y \frac{m(y, \bfz + \bfk - \bfe(x), t-s)}{m(x, \bfz, t-s)}.
\end{equation}  
In other words, compared to the original process, at any time $s \leq t$, transitions along the distinguished lineage are biased in favor of those which lead to a larger $\psi$-average offspring at the final time $t$. 
However, due to the density-dependence of division rates, it is necessary to keep track of the population state $\zeta^{(t)}$. Again, transitions need to be biased, in order to account for the modified behavior of the distinguished individual when compared to the original process. As a consequence, when the population is in state $\bfz$ and the spine is of type $x$ at time $s \leq t$, an individual of type $y$ other than the spine is replaced by offspring $\bfk$ at rate 
\begin{equation}
\label{eq:hrho}
	\widehat{\rho}^{(t)}_\bfk(s, y, x,\bfz) = \tau_\bfk(y, \bfz) \frac{m(x, \bfz + \bfk - \bfe(y), t-s)}{m(x, \bfz, t-s)}.
\end{equation}
Here, the bias favors those transitions which lead to a more favorable environment for the spine, \emph{i.e.} a population composition in which the $\psi$-average of the spine's offspring is high. 

We work under the following assumption, which ensures that the reproduction rates of the spinal process are bounded.
\begin{assumption} 
\label{hyp:bound-m} 
We consider that the following inequality holds: for any $t \geq 0$,
\begin{equation*}
	\max_{x,y \in \X} \; \sup_{\bfz, \bfk \in \Z} \; \sup_{s \in [0,t]} \Big| \frac{m(y, \bfz + \bfk - \bfe(x), s)}{m(x, \bfz, s)} \Big| < \infty.
\end{equation*}
\end{assumption}
In particular, $(Y^{(t)}(s), \zeta^{(t)}(s))_{s \leq t}$ can then be characterized as the unique solution of a stochastic differential equation, as will become apparent in upcoming Section \ref{sec:proof-thm1}.

\begin{remark}
Notice that $m$ can always be bounded from above by the average number of individuals alive in a branching process without deaths, where each individual is replaced by $n > 0$ children at rate $\sum_{k : \|k\|_1 = n}\| \tau_\bfk \|_\infty$. Similarly, $m(x, \bfz, t)$ can be bounded from below by the probability that $u_x(0)$ has not reproduced nor died until time $t$, which is greater than a positive constant independent from $(x, \bfz)$ if for instance $\| \tau_\mathbf{0} \|_\infty < \infty$. This condition thus suffices to ensure that Assumption \ref{hyp:bound-m} holds.
\end{remark}

We are now ready to state our main result. With slight abuse of notation, $\E_{x, \bfz}$ denotes the expectation conditionally on the event $(Y^{(t)}(0), \zeta^{(t)}(0)) = (x, \bfz)$. 
\begin{theorem}
\label{thm:spine}
	Consider that Assumptions \ref{hyp:avg} and \ref{hyp:bound-m} hold. Then for any $t \geq 0$ and any measurable function $F: \D([0,t], \configS) \to \R_+$, for any $ \bfz \in \Z$, 
	{ \small
	\begin{equation}
	\label{eq:spine}
		\popE[\sum_{\substack{u \in \generation{t}}} \psi(x_u(t), Z(t)) F((x_u(s), Z(s))_{s \leq t})] = \sum_{x \in \X} \bfz_x m(x, \bfz, t) \spineE[F((Y^{(t)}(s), \zeta^{(t)}(s))_{s \leq t})].
	\end{equation}
	}
\end{theorem}

This result shows that the reproduction rates defined in Equations \eqref{eq:rho} and \eqref{eq:hrho} yield the appropriate survivorship bias: reproduction events are favored if and only if they increase the $\psi$-average of the offspring at sampling time. In principle, the survivorship bias may thus be computed explicitly, giving access to more precise interpretation. This however requires fine analysis of the function $m$, which at first glance may not seem immediate from Equation \eqref{eq:mpsi}. This motivates the following remark. 

\begin{remark}
\label{rk:ode}
It will become apparent in Lemma \ref{lem:dsmpsi} of Section \ref{sec:proof-thm1} that if 
\[ \exists \alpha > 2 : \quad \sum_{\bfk} \| \bfk \|_1^\alpha \| \tau_\bfk \|_\infty < \infty,\]
then $m$ solves the following linear ODE system, with G defined by Equation \eqref{eq:defG}: 
\begin{equation}
\label{eq:ode-m}
	\begin{cases}
	\partial_t m(x, \bfz,t) = G (m(\cdot,t)) (x,\bfz) \; & \forall (x, \bfz) \in \configS, \forall t > 0, \\
	m(x, \bfz, t = 0) = \psi(x, \bfz) & \forall (x, \bfz) \in \configS. 
	\end{cases}
\end{equation}
If the population size is bounded (\emph{e.g.} epidemiological models neglecting demography; carrying capacity in ecology), the ODE system is finite-dimensional and admits a unique solution. Thus, $m$ can be computed numerically, and may even admit a simple analytical expression as illustrated in Section \ref{sec:toymodel}. Otherwise, the system is finite countable, in which case further assumptions need to be checked to ensure that it admits a unique solution \cite[Chapter 6]{deimlingOrdinaryDifferentialEquations1977}.
\end{remark}

Finally, let us compare the obtained time-inhomogeneous spinal process with the time-homogeneous spine introduced in Section \ref{sec:psi-spine} \cite{bansayeSpineInteractingPopulations2024}.

\begin{remark}
Both the time-homogeneous and time-inhomogeneous spinal constructions are similar in spirit, as the former relies on the $h$-transform, whereas the latter may be regarded as a time-inhomogeneous $m$-transform. Second, in the special case where $\psi$ is an eigenfunction of $G$, a brief computation shows that as expected, Equation \eqref{eq:spine} amounts to the Feynman-Kac formula of \cite[Proposition 1]{bansayeSpineInteractingPopulations2024}. 
\end{remark}

\subsection{Application to a simple density-dependent population process}
\label{sec:toymodel}

Let us end this section by illustrating how the time-inhomogeneous spinal process may be used to gain insight into survivorship bias, in the case where the population size is bounded by a constant $K$, and let 
\[ \configS_K = \{(x, \bfz) \in \X \times \Z: \; \| \bfz \|_1 \leq K, \, \bfz_x \geq 1\}. \]

Generally speaking, in order to make use of the time-inhomogeneous spinal process, the key is the computation of $m$. The set $ \configS_K $ being of finite dimension $d$, Equation \eqref{eq:ode-m} implies that $m$ is characterized as the unique solution of a finite-dimensional linear system of ODEs. More precisely, with slight abuse of notation, let $G$ denote the matrix form of generator $G$ defined by Equation \eqref{eq:defG}. It then follows that for $t \geq 0$, the vector $m(t) = (m(x, \bfz,t) : (x, \bfz) \in \configS_K)$ is given by
\begin{equation*}
m(t) = \mathrm{e}^{Gt} \psi. 
\end{equation*}
While this implies that $m$ can always be computed numerically, there are cases for which an analytical expression of $m$ is achievable. In particular, assume that $G$ is diagonalizable, with linearly independent eigenvectors $v_1, \dots, v_d$ associated to eigenvalues $\lambda_1, \dots, \lambda_d$. Then
\begin{equation*}
	m(t) = \sum_{k=1}^d c_k e^{\lambda_k t} v_k, 
\end{equation*}
where the constants $c_1, \dots, c_d$ are such that $m(0) = \psi$. Throughout the following, we make use of this observation in order to compute the reproduction rates of the time-inhomogeneous spinal process.
\medskip

Let us illustrate the computation of the time-inhomogeneous spinal process on a toy model describing a population of at most two particles, which at each time are either of type $A$ or $B$. Indeed, while computations are achievable for larger populations and state spaces, we restrict ourselves to this minimalistic population in order to keep $\configS_K$ small, so that we can exhibit and comment all transition rates of the time-inhomogeneous spinal process. Our population is thus described by $\bfz = (\bfz_A, \bfz_B)$, where $\bfz_x$ counts the number of particles in state $x$.

The particles behave as follows. If the population corresponds to a single particle of type $A$, the latter may give birth to another particle of the same type at rate $b$. Whenever there are two particles of type $A$, each may die due to competition at rate $a$, or escape competition by switching to state $B$ at rate $c_A$. Finally, a particle of type $B$ can only switch back to state $A$, at rate $c_B$. Throughout the following, we assume that $b=a$ and $c_A = c_B$, for ease of computation.  In other words, the state space of the population with a distinguished particle is described by the ordered set $\{(A, 1, 0), (A, 2, 0), (A,1,1), (B, 1, 1) \}$, and the dynamics are characterized by four reproduction rates:
\begin{equation*}
\tau_{(2,0)}(A, 1, 0) = \tau_{(0,0)}(A, 2,0) = b, \text{ and } \tau_{(0,1)}(A, 2,0) = \tau_{(1,0)}(B, 1,1) = c. \end{equation*}

In this case, the generator $G$ takes the matrix form
\[G = \begin{pmatrix}
-b & 2b & 0 & 0 \\
b & -2(b+c) & c & c \\
0 & c & -c & 0 \\
0 & c & 0 & -c \\
\end{pmatrix}.
 \] 
In particular, $G$ is diagonalizable, with non-positive eigenvalues as the population size is bounded. Hence $m$ can be computed as mentioned previously. For instance, in the case $\psi = \mathbf{1} = (1,1,1,1)$, we obtain that
\[m(t) = \frac{4}{5} \begin{pmatrix} 2 \\ 1 \\ 1 \\ 1 \\ \end{pmatrix} 		
		- \frac{\lambda_-}{10\sqrt{\Delta}} \mathrm{e}^{-\frac{\lambda_+ t}{2}} \begin{pmatrix} -\frac{-3b + 3c - \sqrt{\Delta}}{2c} \\ -\frac{-3b + c + \sqrt{\Delta}}{2c} \\ 1 \\ 1 \\ \end{pmatrix}  		
		+ \frac{\lambda_+}{10\sqrt{\Delta}} \mathrm{e}^{-\frac{\lambda_- t}{2}} \begin{pmatrix} -\frac{-3b + 3c + \sqrt{\Delta}}{2c} \\ -\frac{-3b + c - \sqrt{\Delta}}{2c} \\ 1 \\ 1 \\ \end{pmatrix},\]
with $\Delta = 9b^2 + 9 c^2 - 2bc$ and $\lambda_{\pm} = 3b + 3c \pm \sqrt{\Delta}$.

It thus is possible to explicitly compute the reproduction rates of the time-inhomogeneous spinal process by Equations \eqref{eq:rho} and \eqref{eq:hrho}. For example, assume that sampling occurs at time $t$. If at time $s$, the spine is the sole particle alive and of type $A$, it gives birth to a second particle of type $A$ at rate
\begin{equation*}
\begin{aligned}
\rho^{(t)}_{A,(2,0)}&(s, A, 1,0)= 2b \frac{m(A, 2, 0, t-s)}{m(A, 1, 0, t-s)} \\
&= 2b \frac{16 c \sqrt{\Delta} + \lambda_- (3b + c + \sqrt{\Delta})e^{-\frac{\lambda_+(t-s)}{2}} - \lambda_+(3b + c - \sqrt{\Delta}) e^{-\frac{\lambda_-(t-s)}{2}}}{32 c \sqrt{\Delta} + \lambda_- (-3b + 3c - \sqrt{\Delta})e^{-\frac{\lambda_+(t-s)}{2}} - \lambda_+(-3b + 3c + \sqrt{\Delta}) e^{-\frac{\lambda_-(t-s)}{2}}}.
\end{aligned}
\end{equation*}
The other reproduction rates of the time-inhomogeneous spinal process can be computed analogously. 

\begin{figure}[tb]
   \centering
   \includegraphics[width = 0.8\textwidth]{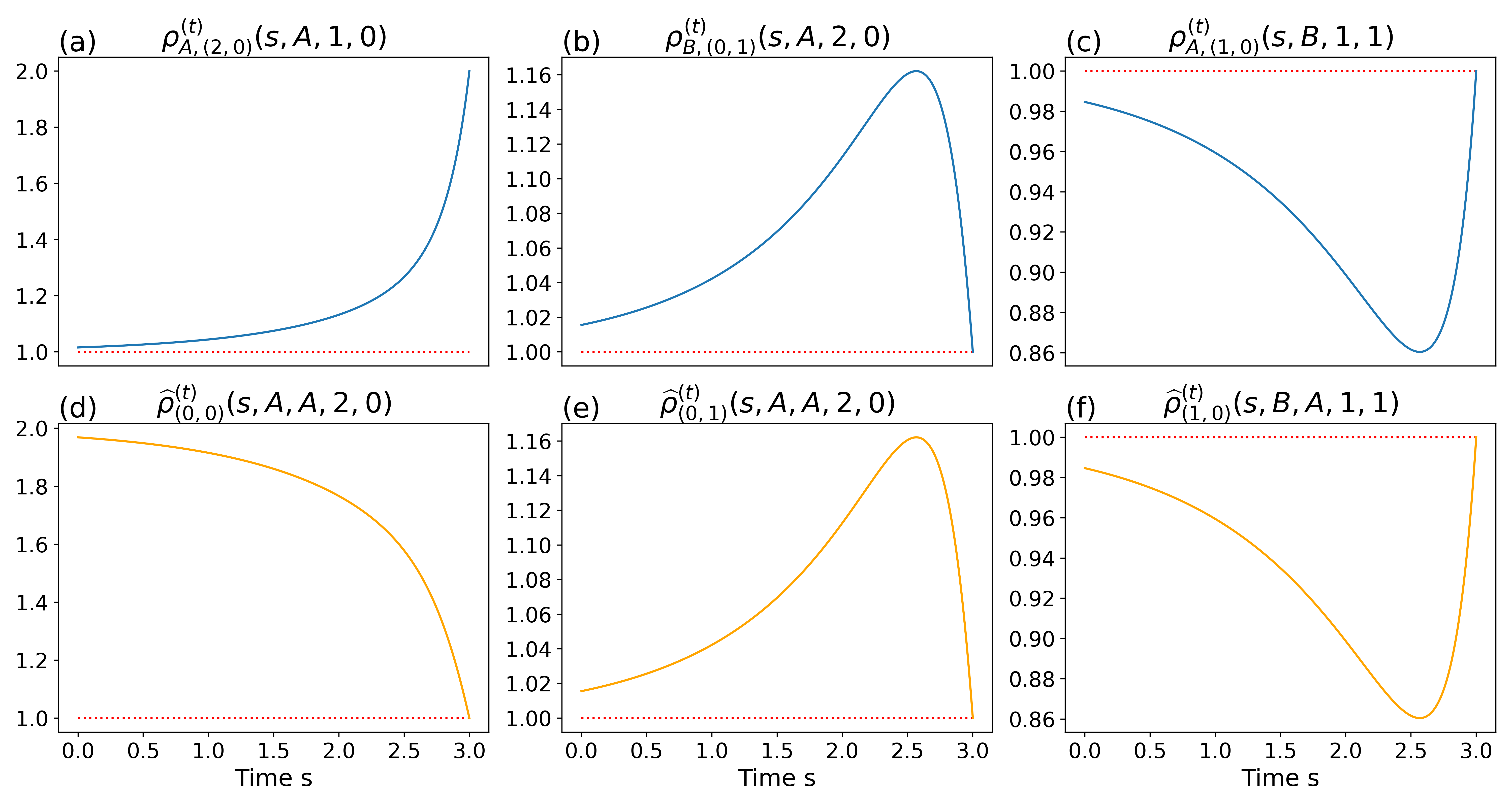} 
   \caption{Time-inhomogeneous reproduction rates of the spine (blue, top) and outside of the spine (orange, bottom). Red dotted lines indicate the unbiased transition rates. Parameters: $t=3$, $b=1$, $c=2$.}
   \label{fig:spine-bd-new}
\end{figure}

Figure \ref{fig:spine-bd-new} illustrates the reproduction rates of the time-inhomogeneous spinal process for the (arbitrary) parameter choice $t=3$, $b=1$ and $c=2$. Recall that in the original process, the sum of the weights $\psi$ of the descendants of a given particle is maximal if the latter has two descendants alive at time $t$, when sampling occurs. Hence, both the birth rate of the spine (panel a) and the death rate of the other particle (panel d) are inflated, as the latter needs to die for the spine to be able to give birth. In addition, whenever there are two particles, transitions from state $A$ to state $B$ are favored (panels b, e), whereas transitions from state $B$ to state $A$ are repressed (panels c, f). In other words, when the population is of size two, the auxiliary process favors the population state $(1,1)$, thus limiting the competitive pressure exerted on the spine.

This illustrates our contributions on capturing the survivorship bias associated to sampling in density-dependent populations of finite size.

\section{The empirical distribution of ancestral lineages in large populations}
\label{sec:res-lln}

In this section, we are interested in empirical distribution of ancestral lineages in large populations, corresponding to the case $\psi = 1$. Indeed, as detailed in Section \ref{sec:pop-lln}, when the order of magnitude of population size $K$ grows to infinity, the population trait distribution converges under classical assumptions to a deterministic limit $z$, characterized as the unique solution of a dynamical system. Building on the large population limit of the time-homogeneous spine construction \cite{bansayeSpineInteractingPopulations2024} which we recall in Section \ref{sec:psi-spine}, we show that it is possible to sample in the limit $z$ with a time-inhomogeneous spine construction, analogously to Theorem \ref{thm:spine}. This allows to explicitly capture the large population survivorship bias. In addition, we quantify the approximation error made by sampling in the large population limit instead of the original finite-size population process.

\subsection{Large population limit of the population process}
\label{sec:pop-lln}

We first focus on the population process, when the initial population size is large and the population process rescaled by the order of magnitude $K$ of the population size. 

Let $D = \mathrm{Card}(\X)$ and consider a family $(\tau_\bfk, \bfk \in \N^{D})$ of functions $\tau_\bfk: \X \times \R_+^D \to \R_+$ satisfying the following hypothesis.
\begin{assumption} 
\label{hyp:lipschitz} 
The following conditions hold:
\begin{enumerate} 
	\item[(i)] The set $J = \{(x,\bfk) \in \X \times \N^{D}: \tau_\bfk(x, \cdot) \neq 0\}$ is finite.
	\item[(ii)] For every $(x,\bfk) \in J$, $\tau_\bfk(x, \cdot)$ is bounded and Lipschitz continuous on $\R_+^D$.
\end{enumerate}
\end{assumption}
For $K \geq 1$, we are interested in the population process $X^K$ where an individual of type $x$ in a population of composition $\bfz \in \Z_K$ is replaced by offspring $\bfk$ at rate $\tau_\bfk(x, \bfz/K)$. 

Introducing the large population limit of the population process requires some notations. 
Define $A(\bfz) = (A_{x,y}(\bfz))_{x,y \in \X}$ for $\bfz \in  \R_+^D$ by $A_{x,y}(\bfz) = \sum_{\bfk} (\bfk_y - \mathbf{e}_y(x)) \tau_\bfk(x, \bfz)$. 
In particular, Assumption \ref{hyp:lipschitz} ensures that $\bfz \mapsto \bfz A(\bfz)$ is locally Lipschitz continuous on $ \R_+^D$. For any $\bfz \in (0, +\infty)^D$, there exists a unique solution $z$ to the differential equation 
\begin{equation}
\label{eq:ode}
	z'(t) = z(t)A(z(t)), \;\; z(0) = \bfz.
\end{equation}
Then the sequence of processes $(Z^K, K \geq 1)$ with initial conditions $Z^K(0) = \lfloor K \bfz \rfloor / K$ converges uniformly in probability to $z$ on finite time intervals \cite[Theorem 3.1, Chapter 11]{ethierMarkovProcessesCharacterization1986}. 

The goal of this Section is to study the empirical distributions of ancestral lineages in this large population limit, which essentially corresponds to sampling in the deterministic limit $z$ itself. As previously, our aim is to explicitly capture the associated survivorship bias. In order to achieve this, we make use of the large population limit of the time-homogeneous spine construction, which we recall below.

\subsection{The time-homogeneous spine construction}
\label{sec:psi-spine}

In this paragraph, we briefly introduce the time-homogeneous spine construction \cite{bansayeSpineInteractingPopulations2024}, as the proof of forthcoming Theorem \ref{thm:lln-spine} builds on it. We first focus on populations of finite size. For $\psi = 1$, individuals other than the spine behave exactly as in the original population process. Whenever the spine is of type $x$ in a population of type distribution $\bfz \in \Z$ such that $\bfz_{x} > 0$, it is replaced by offspring $\bfk$ and becomes of type $y$ at rate $\bfk_y \tau_\bfk(x, \bfz)$. 

The large population limit of the time-homogeneous spinal process is derived as follows. Let $(Y^K,\zeta^K)$ be the spine construction with initial condition $(x, \bfz^K)$ in a population of size $K$. Then \cite[Proposition 7]{bansayeSpineInteractingPopulations2024} ensures that $(Y^K,\zeta^K)_{K \geq 1}$ converges in law, on finite time intervals, to $(\Upsilon, z)$ where $z$ is the unique solution to \eqref{eq:ode}, and $\Upsilon$ is a time-inhomogeneous $\X$-valued Markov jump process which, at time $t$, transitions from $x$ to $y$ at rate 
\begin{equation*}
	\sum_{\bfk} \tau_\bfk(x, z(t))\bfk_y.
\end{equation*} 
In particular, the many-to-one formula associated to sampling in the large population limit relies on an exponential weighting of spinal trajectories, as expected. 
For $(x, \bfz) \in \X \times \Z$, let 
\begin{equation*}
\lambda(x, \bfz) =  \sum_{\bfk : (x, \bfk) \in J} (\| \bfk \|_1 - 1)\tau_\bfk(x, \bfz).
\end{equation*}  
For $t \geq 0$, define
\begin{equation}
\label{eq:tilw}
\mathcal{W}(t) =\exp\left(\int_0^t \lambda(\Upsilon(s), z(s))ds\right).
\end{equation}
With these notations, \cite[Proposition 7]{bansayeSpineInteractingPopulations2024} provides the following many-to-one formula. For any $\bfz \in (0, +\infty)^D$, letting $\bfz^K = \lfloor K \bfz \rfloor / K$,
{ \small
	\begin{equation}
	\label{eq:lln-homogeneous-spine}
	\left| \frac{1}{K} \E_{\bfz^K}\left[\sum_{u \in \generation{t}} F((x_u(s), Z^K(s))_{s \leq t})\right] - \sum_{x \in \X}  \bfz_x \E_{x, \bfz}\left[\mathcal{W}(t)F((\Upsilon(s), z(s))_{s \leq t})\right] \right| \xrightarrow[K \to \infty]{} 0.
	\end{equation}
	}

We are now ready to state our contributions.

\subsection{Survivorship bias in the large population limit}

Les us introduce the large population limit of the time-inhomogeneous spinal process.The general idea is that the impact of the spine on the spinal population becomes negligible in the large population limit, so that in large populations, $\zeta^{(t)}$ is well approximated by $z$. As a consequence, only the dynamics of the spinal individual need to be described. 

Analogously to the finite population case, we expect that the survivorship bias depends on the expected number $\m(x, \bfz, t)$ of individuals alive at time $t \geq 0$, who descend from an original ancestor of type $x$ in a population of type distribution $\bfz \in \R_+^D  \setminus{0}$. In order to capture this quantity, we make use of the large population limit of the time-homogeneous spine \cite{bansayeSpineInteractingPopulations2024} provided by Equation \eqref{eq:lln-homogeneous-spine}. Letting $\overline{\configS} = \X \times (0, +\infty)^D$, we obtain that for any $t \geq 0$ and $(x, \bfz) \in \overline{\configS}$,
\begin{equation*}
 \m(x, \bfz, t) = \spineE[\mathcal{W}(t)],
\end{equation*}
where $\mathcal{W}$ is the stochastic weight associated to the large population limit of the time-homogeneous spine given by Equation \eqref{eq:tilw}. 

We are now ready to define the time-inhomogeneous spinal process. We want to sample at time $t > 0$ in a population of positive initial condition $\bfz \in (0, +\infty)^D$. At time $s \leq t$, the population is of composition $z(s)$, and we let $\Upsilon^{(t)}(s)$ be the type of the spine. Then $\Upsilon^{(t)}$ is a time-inhomogeneous $\X$-valued Markov jump process which, at time $s \in [0,t)$, transitions from $x$ to $y$ at rate 
\begin{equation*}
	\sum_{\bfk} \rho^{(t)}_{\bfk, y}(x, z(s), s),
\end{equation*} 
where 
\begin{equation}
\label{eq:rho-lln}	
	\rho^{(t)}_{\bfk, y}(x, z(s), s) = \tau_\bfk(x, z(s)) \bfk_y \frac{\m(y, z(s), t-s)}{\m(x, z(s), t-s)}.
\end{equation}

Our main contribution lies in upcoming Theorem \ref{thm:lln-spine}, which provides a many-to-one formula for sampling in the large population limit with the time-inhomogeneous spine, and quantifies the associated speed of convergence. Throughout the following, we let $\E_{x, \bfz}$ denote the expectation conditionally on $(\Upsilon^{(t)}(0), z(0)) = (x, \bfz)$. Finally, we introduce the function set
{\small
\begin{equation*}
\begin{aligned}
	\mathfrak{L} = \Big\{F &: \mathbb{D}([0,t], \X)  \times \mathbb{D}([0,t], \R_+^D) \to \R \text{ bounded s.t. } \exists L : \forall x \in \mathbb{D}([0,t], \X), \; \forall \bfz_1, \bfz_2 \in \mathbb{D}([0,t], \R_+^D), \\
	& |F((x(s), \bfz_1(s))_{s \leq t}) - F((x(s), \bfz_2(s))_{s \leq t})| \leq L \sup_{s \in [0,t]} \|\bfz_1(s) - \bfz_2(s)\|_1 \Big\}.
\end{aligned}
\end{equation*}
}
We are now ready to state the main result of this Section, whose proof is postponed to Section \ref{sec:proof-thm2}. 

\begin{theorem}
\label{thm:lln-spine}
	Let $t > 0$ and $F \in \mathfrak{L}$. Under Assumption \ref{hyp:lipschitz}, for any $\bfz \in (0, +\infty)^D$, there exists $C > 0$ such that for any $K \geq 1$, letting $\bfz^K = \lfloor K \bfz \rfloor/K$, 
	{ \small
	\begin{equation*}
	\left| \frac{1}{K} E_{\bfz^K}\left[\sum_{u \in \generation{t}} F((x_u(s), Z^K(s))_{s \leq t})\right] - \sum_{x \in \X} \bfz_{x} \m(x, \bfz, t) \E_{x, \bfz}\left[F((\Upsilon^{(t)}(s), z(s))_{s \leq t})\right] \right| \leq \frac{C}{K^{1/4}}.
	\end{equation*}
	}
\end{theorem}

This theorem shows that the spinal process $(\Upsilon^{(t)}, z)$ indeed captures the typical lineage of a sampled particle in the large population limit $z$. In particular, the survivorship bias thus is correctly captured by Equation \eqref{eq:rho-lln}. Further, if one samples in the large population limit instead of a population whose size is of order $K$, the approximation error is at most of order $K^{-1/4}$. The error thus may exceed the classical approximation error of the population process by its large population limit, which is of order $K^{-1/2}$. This is due to the fact that the approximation of the finite population spinal process by its large population counterpart is associated to two sources of error: the deterministic dynamics $z$ may be different from the type distribution in the finite population spinal process, or the finite population and large population spines may be of different types. It's the combination of both sources of error which leads to the aforementioned speed of convergence, as will become apparent in Section \ref{sec:proof-thm2}.

As in Theorem \ref{thm:spine}, further analyzing the survivorship bias requires a fine understanding of $\m$. In particular, we can show that $\m$ is a solution to a differential equation, analogously to Remark \ref{rk:ode}. Indeed, define the operator $\mathcal{G}$ acting on differentiable functions $f : \overline{\configS} \to \R$ as follows: for any $(x, \bfz) \in \overline{\configS}$, 
\begin{equation*}	
	\mathcal{G} f(x, \bfz) = \sum_\bfk \tau_\bfk(x, \bfz) \angles{\bfk - \bfe(x)}{f(\cdot, \bfz)} + \sum_{\bfk, y} \bfz_y \tau_\bfk(y, \bfz) \angles{\bfk - \bfe(y)}{\nabla_\bfz f(x, \bfz)}.
\end{equation*}
With these notations, we obtain the following result, whose proof is provided in Appendix \ref{appdx:mups}.

\begin{lemma}
\label{lemma:mups}
Under Assumption \ref{hyp:lipschitz}, $\m$ is a solution to the following system of partial differential equations: 
\begin{equation}
\label{eq:pde-m}
\begin{aligned}
	\text{For all } x \in \X \text{, for almost every } t \geq 0 \text{ and } \bfz \in \R_+^D, \\
	\partial_t \m(x,\bfz, t) = \mathcal{G}(\m(\cdot, t))(x, \bfz), \\
	\text{and for any } (x, \bfz) \in \overline{\configS}, \quad \m(x, \bfz, t = 0) = 1.
\end{aligned}
\end{equation}
\end{lemma}
  
In particular, we expect that any weak solution to Equation \eqref{eq:pde-m} is actually a strong solution, and that this solution is unique. Indeed, Equation \eqref{eq:pde-m} defines a hyperbolic system of PDEs. A classical scheme of proof would thus consist to show (1) uniqueness of weak solutions, and (2) existence of a strong solution, which is also a weak solution, see for instance \cite{benzoni-gavageLINEARCAUCHYPROBLEM2006}. Equation \eqref{eq:pde-m} is thus expected to entirely characterize $\m$. This opens the door to further analytical or numerical study of the survivorship bias.

\section{Proofs}
\label{sec:proofs}

\subsection{Proof of Theorem \ref{thm:spine}}
\label{sec:proof-thm1}

The proof proceeds in three steps. First, we show that the time-inhomogeneous spinal process is well defined, and compute the infinitesimal generator associated to its first moment semigroup. Second, we establish the many-to-one formula under the assumption that the population size is bounded. In order to achieve this, we essentially reduce Equation \eqref{eq:spine} to an equality between two time-inhomogeneous semigroups. As the population size is assumed to be bounded, the problem is finite-dimensional and it suffices to show that the semigroups share the same generator to ensure that they are identical. Third, in order to relax this assumption, we consider a localization argument to finally establish Theorem \ref{thm:spine}.

\subsubsection{Existence and uniqueness of the time-inhomogeneous spinal process} We start by characterizing $(Y^{(t)}(s), \zeta^{(t)}(s))_{s \leq t}$ as the unique solution of a stochastic differential equation. In order to do so, we let $Y^{(t)}(s) \in \{\bfe(x): x \in \X\}$ for any $s \in [0,t]$, where $Y^{(t)}(s) = \bfe(x)$ means that the spine is of type $x$. Define $E = \R_+ \times \configS$, and consider two independent Poisson point processes $Q$ and $\widehat{Q}$ on $\R_+ \times E$, of density $dr \otimes d\theta \otimes n(dy, d\bfk)$ where $dr, d\theta$ denote the Lebesgue measure and $n$ the counting measure on $\configS$. Here, we suppose that $Q$ and $\widehat{Q}$ are defined on the same probability space as and independently from $(Y^{(t)}(0), \zeta^{(t)}(0))$, whose law is assumed to be given. Then, for any $s \in [0,t]$, 

{\footnotesize
\begin{equation}
\label{eq:def-sde}
\begin{aligned}
	&Y^{(t)}(s)= Y^{(t)}(0) + \int_0^s \int_E \setind{\theta \leq \rho^{(t)}_{y, \bfk}(r, Y^{(t)}(r-),\zeta^{(t)}(r-))}(\bfe(y) - Y^{(t)}(r-)) Q(dr, d\theta, dy, d\bfk), \\
	&\zeta^{(t)}(s) = \zeta^{(t)}(0) + \int_0^s \int_{E} \setind{\theta \leq \rho^{(t)}_{y, \bfk}(r, Y^{(t)}(r-),\zeta^{(t)}(r-))}(\bfk - Y^{(t)}(r-)) Q(dr, d\theta, dy, d\bfk) \\
	&+ \int_0^s \int_{E} \setind{\theta \leq (\zeta^{(t)}_y(r-) - \setind{Y^{(t)}(r-) = y}) \widehat{\rho}^{(t)}_{\bfk}(r, y,Y^{(t)}(r-),\zeta^{(t)}(r-))}(\bfk - \bfe(y)) \widehat{Q}(dr, d\theta, dy, d\bfk).
\end{aligned}
\end{equation}
}
\begin{remark}
	Throughout the following, in order to simplify notations, we will make no distinction between the sets $\X$ and $\{\bfe(x): x \in \X\}$, based on the natural bijection between the two sets. For example, $Y^{(t)}(s) = x$ is equivalent to $Y^{(t)}(s) = \bfe(x)$. Similarly, to every real-valued function $f$ on $\configS$, we assign a function $\hat{f}$ on $\{\bfe(x): x \in \X\}$ by $\hat{f}(\bfe(x)) = f(x)$, the map $f \mapsto \hat{f}$ being a bijection between the sets of real-valued functions on $\configS$ and on $\{\bfe(x) : x \in \X\}$. Thus, we will always consider $Y^{(t)}$ to take values in the Skorokhod space $\D([0,t],\X)$, unless mentioned otherwise.
\end{remark}

The following result shows that the process $(Y^{(t)}(s), \zeta^{(t)}(s))_{s \leq t}$ is now well defined, and additionally provides its semi-group $R^{(t)} = (R^{(t)}_{r,s}, r \leq s \leq t)$. We recall that the latter is characterized by its action on bounded functions $f$ on $\configS$: for $0 \le r \leq s \leq t$ and $(x,\bfz) \in \configS$,
\begin{equation*}
	R^{(t)}_{r,s}(x, \bfz) = \E[f(Y^{(t)}(s),\zeta^{(t)}(s)) | (Y^{(t)}(r),\zeta^{(t)}(r)) = (x, \bfz)].
\end{equation*}
\begin{proposition}
\label{prop:psi-process}
Under Assumptions \ref{hyp:avg} and \ref{hyp:bound-m}, Equation \eqref{eq:def-sde} admits a unique strong solution $(Y^{(t)}, \zeta^{(t)})$ in the Skorokhod space $\D([0,t], \configS)$. Its semi-group $R^{(t)}$ satisfies: 
\begin{equation}
\label{eq:defR}
	\forall 0 \leq r \leq s\leq t, \quad R^{(t)}_{r,s} = \e^{\int_r^s A^{(t)}_\tau d\tau},
\end{equation}
where the operator $A^{(t)}$ is characterized by its action on non-negative functions $f$ on $\configS$. For any $s \in [0,t]$ and $(x,\bfz) \in \configS$,
\begin{equation}
\label{eq:defA}
\small
	A^{(t)}_s f(x, \bfz) = m(x,\bfz, t-s)^{-1} \left( G(m(\cdot, t-s)f(\cdot))(x,\bfz) - G(m(\cdot,t-s))(x,\bfz) f(x, \bfz)  \right).
\end{equation}
\end{proposition}

The proof is classical, and postponed to Appendix \ref{appdx:existence}.

Before proceeding, we state a technical result which will serve to establish the many-to-one formula for populations of bounded size.
 
\begin{lemma}
\label{lem:dsmpsi}
	In addition to Assumption \ref{hyp:bound-m}, suppose that 
	\begin{equation}
	\label{eq:hyp-m}
		\exists \alpha > 2: \quad \sum_{\bfk} \| \bfk \|_1^\alpha \| \tau_\bfk \|_\infty < \infty. 
	\end{equation}
	Then for any $(x, \bfz) \in \configS$, the function $t \mapsto m(x, \bfz, t)$ is differentiable on $\R_+$, and we have:
	\begin{equation*}
			\partial_t m(x, \bfz,t) = G (m(\cdot,t)) (x,\bfz).
	\end{equation*}
\end{lemma}

Again, the proof follows classical arguments and is postponed to Appendix \ref{appdx:m}.

\subsubsection{Bounded population size}

For now, let us assume that the population size cannot exceed $K$ individuals, \emph{i.e.}
\begin{equation*}
\forall t \geq 0, \quad Z(t) \in \Z_K = \{\bfz \in (\N \cup \{0\})^\X: \Lone{\bfz} \leq K\} \text{ almost surely}.
\end{equation*}
In particular, we thus consider that the following assumption holds.

\begin{assumption} 
\label{hyp:pop-bounded}
There exists $K \in \N$ such that
\begin{equation*}
\tau_\bfk(x,\bfz) = 0 \quad \text{if} \;\; \Lone{\bfz + \bfk - \bfe(x)} > K.
\end{equation*}
\end{assumption}
Since $\Z_K$ is a finite set, this notably implies that Assumptions \ref{hyp:avg} and \ref{hyp:bound-m} are satisfied, as well as Inequality \eqref{eq:hyp-m}. 

Finally, throughout this section, we let 
\begin{equation*}
\configS_K = \{(x, \bfz) \in \X \times \Z_K: \bfz_x \geq 1\}.
\end{equation*} 

We are now ready to turn to the proof of Theorem \ref{thm:spine} for populations of size at most $K$. The proof comprises several steps. Let $t \geq 0$, and start by introducing the time-inhomogeneous semi-group of interest $P^{(t)} = (P^{ (t)}_{r,s}, r \leq s \leq t)$ through its action on functions $f: \configS_K \to \R_+$. For $s \leq t$, $u_x(s)$ is a chosen individual of type $x$ in $\generation{s}$, if it exists. For any $(x, \bfz) \in \configS_K$ and  $0 \leq r \leq s \leq t$, 
\begin{equation}
\label{eq:defp}
	P^{(t)}_{r,s}f(x, \bfz) = m(x, \bfz, t-r)^{-1} \E[\sum_{\substack{u \in \generation{t}\\ u \succeq u_x(r)}}\psi(x_u(t), Z(t)) f(x_u(s), Z(s)) | X(r) = \mathfrak{X}(\bfz)].
\end{equation}

\begin{lemma}
\label{lem:semi-group}
	Under Assumption \ref{hyp:pop-bounded}, $(P^{ (t)}_{r,s}, r \leq s \leq t)$ defines a conservative, time-inhomogeneous semi-group acting on the set of functions $\{f: \configS_K \to \R_+\}$. 
\end{lemma}

\begin{proof}
	The conservativity of $P^{ (t)}$ follows directly from Equation \eqref{eq:defp} applied to $f \equiv 1$, which shows that $P^{ (t)}1 \equiv 1$. 
	
	Let us now turn to the semi-group property. Let $r \leq \tau \leq s \leq t$, and consider $f: \configS_K \to \R_+$ and $(x, \bfz) \in \configS_K$. Throughout the proof, we let $\mathfrak{X}_0 = \mathfrak{X}(\bfz)$.  By definition of the semi-group, 
	\begin{equation}
	\label{eq:sg-prop}
	\begin{aligned}
		P^{ (t)}_{r,s}f(x, \bfz) &= m(x, \bfz,t-r)^{-1} \E[\sum_{\substack{u \in \generation{t}\\ u \succeq u_x(r)}} \psi(x_u(t), Z(t)) f(x_u(s), Z(s)) | X(r) = \mathfrak{X}_0] \\
		&= m(x, \bfz,t-r)^{-1} \E[\sum_{\substack{v \in \generation{\tau}\\ v \succeq u_x(r)}} \sum_{\substack{u \in \generation{t}\\ u \succeq v}} \psi(x_u(t), Z(t)) f(x_u(s), Z(s)) | X(r) = \mathfrak{X}_0] \\
		P^{ (t)}_{r,s}f(x, \bfz) &=  m(x, \bfz,t-r)^{-1} \E[\sum_{\substack{v \in \generation{\tau}\\ v \succeq u_x(r)}} g(x_v(\tau), Z(\tau)) | X(r) = \mathfrak{X}_0],
	\end{aligned}
	\end{equation}
	where we define the function $g: \configS_K \to \R_+$ by 
	\begin{equation*}
		g(x, \bfz) = \E[\sum_{\substack{u \in \generation{t} \\ u \succeq u_x(\tau)}} \psi(x_u(t), Z(t)) f(x_u(s), Z(s))|X(\tau) = \mathfrak{X}_0].
	\end{equation*}	
	Notice that, for any measurable function $G: \D([0,\tau],\configS_K) \to \R_+$,
	{\small
	\begin{equation}
	\label{eq:ttau}
	\begin{aligned}
		&\E[\sum_{\substack{u \in \generation{t}\\ u \succeq u_x(r)}} \psi(x_u(t), Z(t)) m(x_u(\tau), Z(\tau), t-\tau)^{-1} G((x_u(\sigma), Z(\sigma))_{\sigma \leq \tau})) |X(r) = \mathfrak{X}_0] \\
		&= \E[\!\! \sum_{\substack{v \in \generation{\tau}\\ v \succeq u_x(r)}} \!\! \E[\!\! \sum_{\substack{u \in \generation{t}\\ u \succeq v}} \!\! \psi(x_u(t), Z(t)) |X(\tau)] m(x_v(\tau), Z(\tau), t-\tau)^{-1} G((x_v(\sigma), Z(\sigma))_{\sigma \leq \tau}))  |X(r) = \mathfrak{X}_0] \\
		&= \E[\sum_{\substack{v \in \generation{\tau}\\v \succeq u_x(r)}} G((x_v(\sigma), Z(\sigma))_{\sigma \leq \tau})) |X(r) = \mathfrak{X}_0].
	\end{aligned}
	\end{equation}
	}
	Applying this equality to $G((x_v(\sigma), Z(\sigma))_{\sigma \leq \tau}) = g(x(\tau),Z(\tau))$ finally yields the desired semi-group property:
	\begin{equation*}
		P^{ (t)}_{r,s}f(x, \bfz) = P^{ (t)}_{r, \tau} P^{ (t)}_{\tau, s} f(x, \bfz).
	\end{equation*}
	This concludes the proof.
\end{proof}

Let us now compute the generator of $(P^{ (t)}_{r,s}, r \leq s \leq t)$.

\begin{lemma}
\label{lem:gen}
	Let $ t \geq 0$. Under Assumption \ref{hyp:pop-bounded}, the generator of the semi-group $(P^{ (t)}_{r,s}, r \leq s \leq t)$ is $(A^{ (t)}_s, s \leq t)$. 
\end{lemma}

\begin{proof}
	Consider $f: \configS_K \to \R_+$. Let $(x, \bfz) \in \configS_K$ and $t \geq 0$. For any $0 \leq s \leq t$ and $h > 0$ such that $s+h \leq t$, it follows from Equation \eqref{eq:ttau} and the Markov property that 
	\begin{equation*}
		P^{(t)}_{s, s+h}f(x, \bfz) = m(x, \bfz, t-s)^{-1}\popE[\!\!\! \sum_{\substack{u \in \generation{h}\\ u \succeq u_x(0)}} \!\!\! m(x_u(h), Z(h), t-(s+h))f(x_u(h),Z(h))].
	\end{equation*}
	Using Lemma \ref{lem:dsmpsi} as well as the fact that $\configS_K$ is a finite set, we obtain the following Taylor expansion: 
	\begin{equation*}
	\begin{aligned}
	m(x, \bfz, t-s) P^{(t)}_{s, s+h}f(x, \bfz) &= \popE[\!\!\! \sum_{\substack{u \in \generation{h}\\ u \succeq u_x(0)}} \!\!\! m(x_u(h), Z(h), t-s)f(x_u(h),Z(h))] \\
	&+ h \popE[\!\!\! \sum_{\substack{u \in \generation{h}\\ u \succeq u_x(0)}} \!\!\! \partial_s m(x_u(h), Z(h), t-s)f(x_u(h),Z(h))] + o(h).
	\end{aligned}
	\end{equation*}
	As a consequence, 
	\begin{equation*}
	\begin{aligned}
	&m(x, \bfz, t-s) \frac{P^{ (t)}_{s,s+h}f(x, \bfz) - f(x, \bfz)}{h}= \popE[\!\!\! \sum_{\substack{u \in \generation{h}\\ u \succeq u_x(0)}} \!\!\! \partial_s m(x_u(h), Z(h), t-s)f(x_u(h),Z(h))]\\
	& \quad + h^{-1}\Big(\popE[\!\!\! \sum_{\substack{u \in \generation{h}\\ u \succeq u_x(0)}} \!\!\! m(x_u(h), Z(h), t-s)f(x_u(h),Z(h))] - m(x, \bfz,t-s)f(x,\bfz) \Big)
	+ \epsilon(h), 
	\end{aligned}
	\end{equation*}
	where $\epsilon(h)$ is such that $\lim_{h \to 0+} \epsilon(h) = 0$.
	We thus obtain that 
	{ \small
	\begin{equation*}
	\begin{aligned}
		\lim_{h \to 0+} \frac{P^{ (t)}_{s,s+h}f(x, \bfz) - f(x, \bfz)}{h} = m(x, \bfz, t-s)^{-1} \left(G(m(\cdot, t-s)f(\cdot))(x, \bfz) + \partial_s m(x, \bfz, t-s)f(x, \bfz) \right),
	\end{aligned}
	\end{equation*}
	}
	where we recall that $G$ is defined by \eqref{eq:defG}. Lemma \ref{lem:dsmpsi} yields the desired result.
\end{proof}

We finally are ready to establish Theorem \ref{thm:spine} for bounded populations. The proof follows the lines of \cite{marguetUniformSamplingStructured2019}, it is thus only outlined here and we refer to Appendix \ref{appdx:monotone-class} for detail. 

\begin{proposition}
\label{prop:spine-bounded}
Under Assumption \ref{hyp:pop-bounded}, for any $t \geq 0$ and any measurable function $F: \D([0,t], \configS_K) \to \R_+$, for any $ \bfz \in \Z_K$, 
	{ \small
	\begin{equation*}
		\popE[\sum_{\substack{u \in \generation{t}}} \psi(x_u(t), Z(t)) F((x_u(s), Z(s))_{s \leq t})] = \sum_{x \in \X} \bfz_x m(x, \bfz, t) \spineE[F((Y^{(t)}(s), \zeta^{(t)}(s))_{s \leq t})].
	\end{equation*}
	}
\end{proposition}

\begin{proof}
Lemma \ref{lem:gen} implies that the semi-groups $P^{ (t)}$ and $R^{(t)}$ are identical. Using an induction argument, it follows that Equation \eqref{eq:spine} holds for 
\begin{equation*}
F((x(s), \bfz(s))_{s \leq t}) = \prod_{j=1}^k f_j(x(s_j), \bfz(s_j)),
\end{equation*}
where $k \geq 1$, $0 \leq s_1 \leq \dots \leq s_k \leq t$ and $f_1, \dots, f_k: \configS_K \to \R_+$. A monotone class argument finally allows to extend the result to any measurable function $F : \mathbb{D}([0,t], \configS_K) \to \R_+$.
\end{proof}

\subsubsection{Unbounded population size} We are now ready to establish our main result.

\begin{proof}[Proof of Theorem \ref{thm:spine}]

For any $K \in \N$ and $t \geq 0$, we define the following stopping times:
\begin{equation*}
	\sigma_K = \inf \{s \geq 0 : \| Z(s) \|_1 > K \} \quad \text{and} \quad \hat{\sigma}^{(t)}_K = \inf \{s \in [0, t] : \| \zeta^{(t)}(s) \|_1 > K \},
\end{equation*}
with the convention that $\inf \emptyset = +\infty$. In addition, for $t \geq 0$ and $(x, \bfz) \in \configS$, we let 
\begin{equation*}
	m_K(x, \bfz, t) = \popE \Big[ \sum_{\substack{u \in \generation{t} \\ u \succeq u_x(0)}} \psi(x_u(t), Z(t)) \setind{t \le \sigma_K}  \Big]
\end{equation*}
Throughout this proof, for readability, we further make use of the following notations. On the one hand, for $t \geq 0$ and $u \in \generation{t}$, let
\begin{equation*}
 	\couple_u(t) = (x_u(t), Z(t)).
\end{equation*}
Similarly, for $0 \leq s \leq t$, we let
\begin{equation*}
	\couplespine(s) = (Y^{(t)}(s), \zeta^{(t)}(s)).
\end{equation*}

Let $t > 0$, $\bfz \in \Z$ and consider a positive measurable bounded function $F$ on $\D([0,t], \configS)$. With these notations, it follows from Proposition \ref{prop:spine-bounded} that for any $K \geq \| \bfz \|_1$,
\begin{equation*}
	\begin{aligned}
		\popE[\sum_{\substack{u \in \generation{t}}} & \psi(\couple_u(t)) F((\couple_u(s))_{s \leq t}) \setind{t \le \sigma_K}] = \sum_{x \in \X} \bfz_x m_K(x, \bfz, t) \spineE[F((\couplespine(s))_{s \leq t})\setind{t \le \hat{\sigma}^{(t)}_K}].
	\end{aligned}
\end{equation*}
As a consequence, we obtain that 
\begin{equation*}
	\begin{aligned}
	\Big| & \popE[\sum_{\substack{u \in \generation{t}}} \psi(\couple_u(t)) F((\couple_u(s))_{s \leq t})] - \sum_{x \in \X} \bfz_x m(x, \bfz, t) \spineE[F((\couplespine(s))_{s \leq t})] \Big| \\
	& \leq  \popE[\sum_{\substack{u \in \generation{t}}} \psi(\couple_u(t)) F((\couple_u(s))_{s \leq t}) \setind{t > \sigma_K}] \\
	& + \sum_{x \in \X} \bfz_x \Big| m_K(x, \bfz, t) \spineE[F((\couplespine(s))_{s \leq t})\setind{t \le \hat{\sigma}^{(t)}_K}] - m(x, \bfz, t) \spineE[F((\couplespine(s))_{s \leq t})] \Big|.
	\end{aligned}
\end{equation*}

Notice that 
\begin{equation*}
	\popE[\sum_{\substack{u \in \generation{t}}} \psi(\couple_u(t)) F((\couple_u(s))_{s \leq t}) \setind{t > \sigma_K}] \leq \| \psi \|_\infty \|F\|_\infty \popE[\text{Card}(\generation{t}) \setind{t > \sigma_K}].
\end{equation*}
We can couple $Z$ with a branching process $V$ starting from $V(0) = \| \bfz \|_1$ individuals, and within which each individual is replaced by $n > 0$ children at rate $\sum_{\bfk \in \Z: \| \bfk \|_1 = n} \| \tau_\bfk \|_\infty$,
such that $\| Z(s)\|_1 \leq V(s)$ almost surely for any $s \geq 0$. In particular, 
\[  \sigma_K \geq \inf \{ s > 0 : V(s) \geq K \}. \]
Assumption \ref{hyp:avg} ensures non-explosion of $V$. As a consequence, 
\[ \lim_{K \to \infty} \text{Card}(\generation{t}) \setind{t > \sigma_K} = 0 \quad \text{almost surely}\]
and further $\E[\text{Card}(\generation{t})] < \infty$. We thus conclude by dominated convergence that 
\begin{equation}
\label{eq:dc-1}
	\lim_{K \to \infty} \popE[\text{Card}(\generation{t}) \setind{t > \sigma_K}] = 0.
\end{equation}

Next, fix $x \in \X$. We have
\begin{equation*}
\begin{aligned}
	\Big| m_K(x, \bfz, t) \spineE[F&((\couplespine(s))_{s \leq t})\setind{t \leq \hat{\sigma}^{(t)}_K}] - m(x, \bfz, t) \spineE[F((\couplespine(s))_{s \leq t})] \Big| \\
	& \leq m(x, \bfz, t) \|F\|_\infty \P_{x, \bfz}(t > \hat{\sigma}^{(t)}_K) + \|F\|_\infty |m(x, \bfz, t) - m_K(x, \bfz, t)|.
\end{aligned}
\end{equation*}

Proceeding as for Equation \eqref{eq:dc-1}, we obtain that 
\begin{equation*}
	\lim_{K \to \infty} |m(x, \bfz, t) - m_K(x, \bfz, t)| = 0.
\end{equation*}

Similarly, we may couple $\zeta^{(t)}$ with a branching process $W$ starting from $\| \bfz\|_1$ individuals, within which each individual is replaced by $n > 0$ children at rate 
\[\sum_{k \in \Z : \|k\|_1 = n} \max \Big( \sum_y \| \rho^{(t)}_{y, \bfk} \|_\infty, \| \hat{\rho}^{(t)}_{\bfk} \|_\infty \Big)  ,\]
such that $\| \zeta^{(t)}(s)\|_1 \leq W(s)$ almost surely for any $s \in [0,t]$. Assumptions \ref{hyp:avg} and \ref{hyp:bound-m} ensure non-explosion of $W$ on $[0,t]$, which implies that
\begin{equation*}
	\lim_{K \to \infty} \P_{x, \bfz}(t > \hat{\sigma}^{(t)}_K) \to 0.
\end{equation*}

As a consequence, we thus have established that Equation \eqref{eq:spine} holds for any test function $F$ which is positive, measurable and bounded. The result extends to any positive measurable test function $F$ by approximating it with an increasing sequence of simple functions, by monotone convergence.

\end{proof}

\subsection{Proof of Theorem \ref{thm:lln-spine}}
\label{sec:proof-thm2}

Let us start by outlining the proof, which proceeds in three steps. First, we show that $(\Upsilon^{(t)}, z)$ is well defined and characterize its first moment semi-group which is defined by its action on functions $f \in \mathrm{L}^\infty(\overline{\configS})$: for $r \leq s \leq t$ and $(x,\bfz) \in \overline{\configS}$,
\begin{equation*}
	\calR^{(t)}_{r,s}f(x, \bfz) = \E[f(\Upsilon^{(t)}(s),z(s)) | (\Upsilon^{(t)}(r),z(r)) = (x, \bfz)].
\end{equation*}

\begin{proposition}
\label{prop:lln-spine}
Under Assumption \ref{hyp:lipschitz}, $(\Upsilon^{(t)}, z)$ is characterized as the unique strong solution in $\mathbb{D}([0,t], \X) \times \C^1([0,t], \Z)$ to a stochastic differential equation. In addition, $\calR^{(t)}$ is a time-inhomogeneous semi-group of bounded linear operators on $\mathrm{L}^\infty(\overline{\configS})$, whose generator $\A^{(t)}$ is characterized by its action on functions $f \in \mathcal{C}^1(\overline{\configS})$ as follows. For $s \leq t$ and $(x, \bfz) \in \overline{\configS}$, 
\begin{equation*}
\begin{aligned}
	\A^{(t)}_s f(x,\bfz) &= \sum_{\substack{\bfk \in \N^D \\ y \in \X}} \Big( \rho^{(t)}_{\bfk, y}(x, \bfz, s) \big( f(y, \bfz)-f(x,\bfz) \big) + \bfz_y \tau_\bfk(y, \bfz) \angles{\nabla_\bfz f(x, \bfz)}{\bfk - \bfe(y)} \Big).
\end{aligned}
\end{equation*}
\end{proposition}
The proof is classical, and we refer to Appendix \ref{proof:prop-lln-spine} for detail.

Second, we show that $(\Upsilon^{(t)}, z)$ indeed provides a many-to-one formula for sampling in the large population limit. Analogously to the proof of Theorem \ref{thm:spine}, we identify its first moment semigroup with an appropriate time-inhomogeneous semi-group derived from the large population limit of the time-homogeneous spine construction. This allows us to establish the following proposition, whose proof is postponed to Section \ref{sec:proof-prop4}.

\begin{proposition}
\label{prop:link-spines}
	Under Assumption \ref{hyp:lipschitz}, for any measurable function $F: \D([0,t], \X) \times \C^1([0,t], \R_+^D) \to \R_+$, for any $ \bfz \in (0, +\infty)^D$ and $t \geq 0$, 
	{ \small
	\begin{equation}
	\label{eq:link-spines}
		\spineE[\mathcal{W}(t) F((\Upsilon(s), z(s))_{s \leq t})] = \m(x, \bfz, t) \spineE[F((\Upsilon ^{(t)}(s), z(s))_{s \leq t})].
	\end{equation}
	}
\end{proposition}
In particular, letting $\bfz^K = \lfloor K \bfz \rfloor / K$, Equation \eqref{eq:lln-homogeneous-spine} thus implies 
{\small
	\begin{equation*}
	\left| \frac{1}{K} \E_{\bfz^K}\left[\sum_{u \in \generation{t}} F((x_u(s), Z^K(s))_{s \leq t})\right] - \sum_{x \in \X} \bfz_x \m(x, \bfz, t) \spineE[F((\Upsilon ^{(t)}(s), z(s))_{s \leq t})] \right| \xrightarrow[K \to \infty]{} 0.
	\end{equation*}
	}

Third, it remains to quantify the speed of convergence. This is achieved in the following result. 
\begin{proposition}
\label{prop:lln-sample}
	Let $t > 0$ and $F \in \mathfrak{L}$. Fix $\bfz \in (0, +\infty)^D$ and let $\bfz^K = \lfloor K \bfz \rfloor / K$. Under assumption \ref{hyp:lipschitz}, there exists $C > 0$ such that for any sequence $(\varepsilon_K)_{K \geq 0}$ of positive real numbers, for any $K \geq 1$, letting $\delta_K = C(\varepsilon_K + K^{-1/2}\varepsilon_K^{-1})$,
	{ \small
	\begin{equation*}
	\left| \E_{\bfz^K}\left[\sum_{u \in \generation{t}} F((x_u(s), Z^K(s))_{s \leq t})\right] - \sum_{x \in \X} \bfz_{x} \E_{x, \bfz}\left[\mathcal{W}(t) F((\Upsilon(s), z(s))_{s \leq t})\right] \right| \leq \delta_K.
	\end{equation*}
	}
\end{proposition}
The proof proceeds by coupling the time-homogeneous spine in a population of size of ordre $K$ to its large population limit. We refer to Section \ref{sec:proof-prop5} for detail.

Combining Proposition \ref{prop:link-spines} and \ref{prop:lln-sample} for $\varepsilon_K = K^{-1/4}$ finally yields Proposition \ref{prop:lln-spine}.

The rest of this section is structured as follows. We start by establishing a key lemma, which will be useful throughout the proofs of both Propositions \ref{prop:link-spines} and \ref{prop:lln-sample}. We then turn towards the proof of Proposition \ref{prop:link-spines}, and finally conclude by the proof of Proposition \ref{prop:lln-sample}.

\subsubsection{A key lemma}

In this section, we establish a key lemma which states that controlling the difference in population composition is actually sufficient to control the probability that the spinal processes grow apart, in final time. This is achieved by a coupling argument.

Let $\zeta_1$ and $\zeta_2$ be two population processes in $\D([0,t], \R_+^D)$ defined on the same probability space $(\Omega, \mathcal{F}, \P)$. 
Consider a family of independent Poisson Point Processes $(Q_{j}, j \in J^*)$ of intensity the Lebesgue measure on $\R^2_+$, which is also defined on $\Omega$, and for $y_0 \in \X$ let 
\begin{equation*}
\begin{aligned}
	Y_1(t) &= \bfe(y_0) + \sum_{(x, \bfk, y) \in J^*} \int_0^t \int_0^{+\infty} \setind{Y_1(s-)=x, \theta \leq \bfk_y \tau_\bfk(x, \zeta_1(s-))}(\bfe(y) - \bfe(x)) Q_{x, \bfk, y}(ds, d\theta), \\
	Y_2(t) &= \bfe(y_0) + \sum_{(x, \bfk, y) \in J^*} \int_0^t \int_0^{+\infty} \setind{Y_2(s-)=x, \theta \leq \bfk_y \tau_\bfk(x, \zeta_2(s-))}(\bfe(y) - \bfe(x)) Q_{x, \bfk, y}(ds, d\theta). \\
\end{aligned}
\end{equation*}
Notice that $Y_j(t)$ corresponds to the dynamics of the spine, given that the population composition is provided by $\zeta_j$, for $j \in \{1,2\}$. For instance, if $\zeta_1 = \phi(\cdot, \bfz_1)$, we have $(Y_1, \zeta_1)  = (\Upsilon, z)$ with initial condition $ (y_0, \bfz_1)$. 

\begin{remark}
Importantly, we do not assume that the family $(Q_{j}, j \in J^*)$ is independent from $\zeta_1$ and $\zeta_2$. For instance, given a family of independent Poisson Point Processes $(Q_{j}, j \in J \cup J^*)$ of intensity the Lebesgue measure on $\R^2_+$, we may consider $(Y_1, \zeta_1)$ to satisfy upcoming Equation \eqref{eq:psi-spine}, hence recovering the time-homogeneous spinal process $(Y^K, \zeta^K)$ defined in Section \ref{sec:psi-spine}.
\end{remark}

We are now ready to state the main result of this section.

\begin{lemma}
\label{lemma:error}
Let $t > 0$, and assume that there exist two positive sequences $(\varepsilon_K)$ and $(\alpha_K)$  such that for every $K \geq 1$, 
	\begin{equation}
	\label{eq:hyp-zeta}
 	\P\left(\sup_{s \in [0,t]} \|\zeta_1(s) - \zeta_2(s) \|_1 \geq \varepsilon_K \right) \leq \alpha_K. 
	\end{equation}
Under Assumption \ref{hyp:lipschitz}, there exists a constant $C(t) > 0$
	such that for every $K \geq 1$,
	\begin{equation*}
	\P(\forall s \in [0,t], Y_1(s) = Y_2(s)) \geq 1 - C(t)(\alpha_K + \varepsilon_K).
	\end{equation*}
\end{lemma}

\begin{proof}

We are interested in the first instant $T_K$ at which $Y_1$ differs from $Y_2$:
\begin{equation*}
	T_K = \inf\{t \geq 0: Y_1(t) \neq Y_2(t)\}.
\end{equation*}

Let $K \geq 1$ and $\bfz \in (0, +\infty)^D$. For two sets $A$ and $B$, $A \Delta B$ denotes their symmetric difference. For $(x, \bfk, y) \in J^*$, $(x', \bfz) \in \overline{\mathcal{S}}$ and $\theta > 0$, define the event 
\begin{equation*}
	\mathcal{E}_{x, \bfk, y}(\theta, x', \bfz) = \{x' = x, \theta \leq \bfk_y \tau_\bfk(x, \bfz).\}
\end{equation*}

With these notations, we may notice that by the coupling of $Y_1$ and $Y_2$,
{ \small
\begin{equation*}
\begin{aligned}
	\{T_K > t\} & \supseteq \left\{ T_K > t,  \sum_{(x, \bfk, y) \in J^*} \int_0^{t} \int_0^{+\infty} \setind{\mathcal{E}_{x, \bfk, y}(\theta, Y_1(s-), \zeta_1(s-)) \Delta  \mathcal{E}_{x, \bfk, y}(\theta, Y_2(s-), \zeta_2(s-)) } Q_{x, \bfk, y}(ds, d\theta) = 0\right\} \\
	&\supseteq \left\{ T_K > t, \sum_{(x, \bfk, y) \in J^*} \int_0^{t} \int_0^{+\infty} \setind{\mathcal{E}_{x, \bfk, y}(\theta, Y_1(s-),\zeta_1(s-)) \Delta  \mathcal{E}_{x, \bfk, y}(\theta, Y_1(s-), \zeta_2(s-)) } Q_{x, \bfk, y}(ds, d\theta) = 0\right\} \\
	\{T_K > t\} &\supseteq \left\{\sum_{(x, \bfk, y) \in J^*} \int_0^{t} \int_0^{+\infty} \setind{\mathcal{E}_{x, \bfk, y}(\theta, Y_1(s-), \zeta_1(s-)) \Delta  \mathcal{E}_{x, \bfk, y}(\theta, Y_1(s-), \zeta_2(s-)) } Q_{x, \bfk, y}(ds, d\theta) = 0\right\}.
\end{aligned}
\end{equation*}
}

Let $\mathcal{D} = \{(y, \bfk) : \exists x \in \X \text{ s.t. } (x, \bfk) \in J, \bfk_y > 0\}$. Recall from Assumption \ref{hyp:lipschitz} that for any $(x, \bfk) \in J$, $\tau_\bfk(x, \cdot)$ is $L_{x,\bfk}$-Lipschitz continuous. Let $L = \max_{(x, \bfk) \in J} \|k\|_1 L_{x, \bfk}$.
We introduce the event 
\begin{equation*}
	A_K = \left\{\max_{(y, \bfk) \in \mathcal{D}} \sup_{s \in [0,t]}  \bfk_y |\tau_{y, \bfk}(Y_1(s), \zeta_1(s)) - \tau_{y, \bfk}(Y_1(s), \zeta_2(s)) | < L \varepsilon_K \right\}.
\end{equation*}
It follows that
{ \small
\begin{equation}
\label{eq:TK-decomposed}
\begin{aligned}
	\P&(T_K \leq t) = \P(\sum_{(x, \bfk, y) \in J^*} \int_0^{t} \int_0^{+\infty} \setind{\mathcal{E}_{x, \bfk, y}(\theta, Y_1(s-), \zeta_1(s-)) \Delta  \mathcal{E}_{x, \bfk, y}(\theta, Y_1(s-), \zeta_2(s-)) } Q_{x, \bfk, y}(ds, d\theta) \geq 1)\\
	\leq & \; \P(A_K, \sum_{(x, \bfk, y) \in J^*} \int_0^{t} \int_0^{+\infty} \setind{\mathcal{E}_{x, \bfk, y}(\theta, Y_1(s-), \zeta_1(s-)) \Delta  \mathcal{E}_{x, \bfk, y}(\theta, Y_1(s-), \zeta_2(s-)) } Q_{x, \bfk, y}(ds, d\theta) \geq 1) + \P(A_K^C).
\end{aligned}
\end{equation}
}

First, we may notice that Assumption \ref{hyp:lipschitz} ensures that
\begin{equation*}
\begin{aligned}
	A_K^C \subseteq \{\sup_{s \in [0,t]} \|\zeta_1(s) - \zeta_2(s)\|_1 \geq \varepsilon_K \},
\end{aligned}
\end{equation*}
from which we deduce by Equation \eqref{eq:hyp-zeta} that 
\begin{equation}
\label{eq:p1}
	\P(A_K^C) \leq \alpha_K.
\end{equation}

Second, on the event $A_K$, it holds that for any $\theta \geq 0$ and $(x, \bfk, y) \in J^*$,
\begin{equation*}
\begin{aligned}
\mathcal{E}_{x, \bfk, y}&(\theta, Y_1(s-), \zeta_1(s-)) \Delta  \mathcal{E}_{x, \bfk, y}(\theta, Y_1(s-), \zeta_2(s-)) \\
& \subseteq \{\theta \in [\bfk_y \tau_\bfk(Y_1(s-), \zeta_1(s-)) - L\varepsilon_K, \bfk_y \tau_\bfk(Y_1(s-), \zeta_1(s-)) + L\varepsilon_K] \}. 
\end{aligned}
\end{equation*}
As a consequence, 
{\small
\begin{equation*}
\begin{aligned}
	\Big\{ A_K, & \sum_{(x, \bfk, y) \in J^*} \int_0^{t} \int_0^{+\infty} \setind{\mathcal{E}_{x, \bfk, y}(\theta, Y_1(s-), \zeta_1(s-)) \Delta  \mathcal{E}_{x, \bfk, y}(\theta, Y_1(s-), \zeta_2(s-)) } Q_{x, \bfk, y}(ds, d\theta) \geq 1 \Big\} \\
	& \subseteq \Big\{ \sum_{(x, \bfk, y) \in J^*} \int_0^{t} \int_0^{+\infty} \setind{ \theta \in [\bfk_y \tau_\bfk(Y_1(s-), \zeta_1(s-)) - L\varepsilon_K, \bfk_y \tau_\bfk(Y_1(s-), \zeta_1(s-)) + L\varepsilon_K] } Q_{x, \bfk, y}(ds, d\theta) \geq 1 \Big\} . 
\end{aligned}
\end{equation*}
}
Hence, Markov's inequality leads to
\begin{equation}
\label{eq:p2}
\begin{aligned}
	&\P(A_K, \sum_{(x, \bfk, y) \in J^*} \int_0^{t} \int_0^{+\infty} \setind{\mathcal{E}_{x, \bfk, y}(\theta, Y_1(s-), \zeta_1(s-)) \Delta  \mathcal{E}_{x, \bfk, y}(\theta, Y_1(s-), \zeta_1(s)) } Q_{x, \bfk, y}(ds, d\theta) \geq 1) \\
	& \leq \E[ \sum_{(x, \bfk, y) \in J^*} \int_0^{t} \int_0^{+\infty} \setind{ \theta \in [\bfk_y \tau_\bfk(Y_1(s-), \zeta_1(s-)) - L\varepsilon_K, \bfk_y \tau_\bfk(Y_1(s-), \zeta_1(s-)) + L\varepsilon_K] } Q_{x, \bfk, y}(ds, d\theta)] \\
	& \leq Ct \varepsilon_K,
\end{aligned}
\end{equation}
with $C = 2L\mathrm{Card}(J^*)$.
Injecting Inequalities \eqref{eq:p1} and \eqref{eq:p2} into Equation \eqref{eq:TK-decomposed} concludes.
\end{proof}

\subsubsection{Proof of Proposition \ref{prop:link-spines}} 
\label{sec:proof-prop4}

Consider the initial condition $\bfz \in (0, +\infty)^D$ to be fixed. As above, we write $\phi(t, \bfz)$ for the solution to Equation \eqref{eq:ode} evaluated at time $t$. Further, $\Upsilon_\bfz$ denotes the large population limit of the time-homogeneous spine, given that $z(0) = \bfz$. We are interested in the following time-inhomogeneous semi-group, defined through its action on functions $f : \X \to \R_+$. For any $x \in \X$ and $0 \le r \le s \le t$,
\begin{equation*}
	\calP^{(t, \bfz)}_{r, s} f(x) = \m(x, \phi(r, \bfz), t-r)^{-1} \E \Big[\exp\Big(\int_r^t \lambda(\Upsilon_\bfz(u), \phi(u, \bfz)) du \Big) f(\Upsilon_\bfz(s)) | \Upsilon_\bfz(r) = x \Big].
\end{equation*}
 
\begin{proposition}
\label{prop:pt}
	Under Assumption \ref{hyp:lipschitz}, $P^{(t, \bfz)}$ is a time-inhomogeneous, conservative semi-group of bounded linear operators on $L^\infty(\X)$, whose generator $\A^{(t, \bfz)}$ is defined by its action on functions $f : \X \to \R_+$ as follows. For any $x \in \X$ and $0 \le s \le t$,
	\[ \A^{(t, \bfz)} f(x) = \sum_{\bfk \in \N^D} \sum_{y \in \X} \rho^{(t)}_{\bfk, y}(x, \phi(\bfz, s), s) \big( f(y) - f(x) \big).  \]
\end{proposition}

\begin{proof}
	Let $\bfz \in \R_+^D$ be fixed. From the definition of the time-homogeneous large population spine, it follows that 
	\begin{equation}
	\label{eq:aux-m}
	\begin{aligned}
		m(x, \phi(r, \bfz), t-r) &= \E\Big[ \exp \Big(\int_0^{t-r} \lambda(\Upsilon_{\phi(r, \bfz)}(s), \phi(r + s, \bfz))ds \Big) \Big| \Upsilon_{\phi(r, \bfz)}(0) = x \Big] \\
		&= \E\Big[ \exp \Big(\int_r^{t} \lambda(\Upsilon_{\bfz}(s), \phi(s, \bfz))ds \Big) \Big| \Upsilon_\bfz(r) = x \Big].
	\end{aligned}
	\end{equation}
	It thus is clear that $\calP^{(t, \bfz)}_{r,s} 1(x) = 1$ for any $0 \leq r \leq s \leq t$ and $x \in \X$. 
	
	Let us establish the semi-group property, i.e. we want to show that for any $0 \leq r \leq \tau \leq s \leq t$, 
\begin{equation*}
	\calP^{(t, \bfz)}_{r,\tau}\calP^{(t,\bfz)}_{\tau,s} = \calP^{(t, \bfz)}_{r,s}.
\end{equation*}

Throughout the proof, for $s \geq 0$ we write $\Ybb_\bfz(s) = (\Upsilon_\bfz(s), \phi(s, \bfz))$. Further, for $t_1 \leq t_2$ define
\[ \mathcal{W}_\bfz(t_1, t_2) = \exp\Big( \int_{t_1}^{t_2} \lambda(\Ybb(u))du \Big). \]
	
	By definition, we have for $x \in \X$ and $f \in L^{\infty}(\X)$:
\begin{equation*}
\begin{aligned}
m(x, \phi(r, \bfz),& t-r) \calP^{(t, \bfz)}_{r,\tau}\calP^{(t,\bfz)}_{\tau,s} f(x) = \E\Big[ \mathcal{W}_\bfz(r,t) \calP^{(t,\bfz)}_{\tau,s} f(\Upsilon_\bfz(\tau)) \Big| \Upsilon_\bfz(r) = x \Big]. \\
& = \E\Big[ \mathcal{W}_\bfz(r,t) m(\Ybb_\bfz(\tau), t-\tau)^{-1} \E[\mathcal{W}_\bfz(\tau,t) f(\Ybb_\bfz(s)) | \Upsilon_\bfz(\tau)] \Big| \Upsilon_\bfz(r) = x \Big].
\end{aligned}
\end{equation*}
Notice that, almost surely,
\[ \mathcal{W}_\bfz(r,t) = \mathcal{W}_\bfz(r,\tau)\mathcal{W}_\bfz(\tau,t).\]

Let $(\mathcal{F}_s, s \geq 0)$ be the natural filtration associated to $\Ybb_\bfz$. As $\mathcal{W}_\bfz(r,\tau)$ is $\mathcal{F}_\tau$-measurable, we obtain that
\begin{equation*}
\begin{aligned}
m(x, \phi(r, \bfz),& t-r) \calP^{(t, \bfz)}_{r,\tau}\calP^{(t,\bfz)}_{\tau,s} f(x) = \\
&  \E\Big[ \mathcal{W}_\bfz(\tau,t) m(\Ybb_\bfz(\tau), t-\tau)^{-1} \E[\mathcal{W}_\bfz(r,t) f(\Ybb_\bfz(s)) | \Upsilon_\bfz(\tau)] \Big| \Upsilon_\bfz(r) = x \Big].
\end{aligned}
\end{equation*}
It follows from Equation \eqref{eq:aux-m} that $\E[\mathcal{W}_\bfz(\tau,t) | \Upsilon_\bfz(\tau)] = m(\Ybb_\bfz(\tau), t-\tau)$. Since $\mathcal{F}_r \subseteq \mathcal{F}_\tau$, this finally leads to 
\[ m(x, \phi(r, \bfz), t-r) \calP^{(t, \bfz)}_{r,\tau}\calP^{(t,\bfz)}_{\tau,s} f(x) = \E \Big[ \mathcal{W}_\bfz(r, t)f(\Ybb_\bfz(s))  \Big| \Upsilon_\bfz(r) = x \Big], \]
as desired. 

Finally, proceeding as in the proof of Proposition \ref{prop:lln-spine} leads to 
\[ \lim_{h \to 0+} \frac{\calP^{(t, \bfz)}_{s, s+h} f(x) - f(x)}{h} = \A^{(t,\bfz)}f(x). \]
This concludes the proof.
\end{proof}

We are now ready to establish the main result of this section.
\begin{proof}[Proof of Proposition \ref{prop:link-spines}]
	Analogously to previous notations, $\Upsilon^{(t)}_\bfz$ is the time-inhomogeneous large population spine, given that $z(0) = \bfz$. In addition, $\R^{(t, \bfz)}$ is the first moment semi-group associated to the marginal law of $\Upsilon^{(t)}_\bfz$, \emph{i.e.} for any $f : \X \to \R_+$, $x \in \X$ and $0 \le r \le s \le t$,
\begin{equation*}
	\calR^{(t, \bfz)}_{r,s}f(x) = \E[f(\Upsilon^{(t)}_\bfz(s)) | \Upsilon^{(t)}_\bfz(r) = x].
\end{equation*}

	Propositions \ref{prop:lln-spine} and \ref{prop:pt} ensure that the semi-groups $\calP^{(t, \bfz)}$ and $\calR^{(t, \bfz)}$ are identical. Hence Equation \eqref{eq:link-spines} holds for 
\begin{equation*} 
\begin{aligned}
	F : &\;  \D([0,t], \X) \times \C^1([0,t], \Z) \to \R_+ \\
	& (y(s), z(s))_{s\leq t}  \mapsto f(y(s), z(s)),\end{aligned}
\end{equation*}
with $s \in [0,t]$ and $f \in \mathrm{L}^{\infty}(\overline{\configS})$ given. As in the Proof of Theorem \ref{thm:spine}, this suffices to conclude by induction and using a monotone class argument.
\end{proof}

\subsubsection{Proof of Proposition \ref{prop:lln-sample}}
\label{sec:proof-prop5}

It remains to quantify the speed of convergence in Theorem \ref{thm:lln-spine}. Notice that according to Proposition \ref{prop:link-spines}, it is identical to the speed of convergence in Equation \eqref{eq:lln-homogeneous-spine}. By coupling the time-homogeneous spine in a population of size $K$ to its large population limit, we conclude.

\paragraph{Preliminary results on the time-homogeneous spine.} We start by studying the time-homogeneous spinal process in a population of size of order $K$. 

For convenience, we introduce a slight change in the type space allowing us to derive an equation for the spinal population type distribution $\zeta^K$ which does not depend on $Y^K$. More precisely, the type space now becomes $\X^* = \{0,1\} \times \X$. An individual of type $(0,x) \in \X^*$ corresponds to an individual of type $x$ which is not the spine, whereas an individual is of type $(1,x) \in \X^*$ if it is the spine and of type $x$. We let $\zeta^K = (\zeta^K_{i,x}, (i,x) \in \X^*)$ denote the corresponding type distribution of the spinal population. Finally, we define $\proj{\zeta^K} = ((\proj{\zeta^K})_x, x \in \X)$ by 
\begin{equation*}
\left(\proj{\zeta^K}\right)_x = \zeta^K_{1,x} + \zeta^K_{0,x} \quad \forall x \in \X.
\end{equation*}
Throughout the following, with slight abuse of notation, we write $\tau_{\bfk}(x, \bfz)$ instead of $\tau_{\bfk}(x, \proj{\bfz})$ for clarity.

For $(x, \bfk) \in J$ and $(x, \bfk, y) \in J^* = \{(x, \bfk, y) : (x, \bfk) \in J, \bfk_y > 0\}$, we let 
\begin{equation*}
\begin{aligned}
	h_0(x, \bfk) &= \sum_{y \in \X} (\bfk_y - \setind{y=x}) \bfe(0,y) \\
	\text{and} \quad h_1(x, \bfk, y) &= \sum_{w \in \X} (\bfk_w - \setind{y=w}) \bfe(0,w) + \bfe(1,y) - \bfe(1,x).
\end{aligned}
\end{equation*}

As previously, we fix $\bfz \in \R_+^D \setminus \{0\}$, and recall that $\bfz^K = \lfloor K \bfz \rfloor/ K$. For any $x \in \X$, the initial condition $(x, \bfz^K)$ becomes $(x, \overline{\bfz}^K)$ with 
\begin{equation*}
\overline{\bfz}^K = \frac{1}{K} \bfe(1,x) + \sum_{y \in \X} \left((\bfz^K)_{y} - \frac{1}{K} \setind{x=y} \right) \bfe(0,y). 
\end{equation*}
Given a family of independent Poisson Point Processes $(Q_{j}, j \in J \cup J^*)$ of intensity the Lebesgue measure on $\R^2_+$, the process $(Y^K, \zeta^K)$ can then be defined as follows:
{ \small
	\begin{equation}
	\label{eq:psi-spine}
	\begin{aligned}
		Y^K(t) &= \bfe(y_0) + \sum_{(x, \bfk, y) \in J^*} \int_0^t \int_0^{+\infty} \setind{Y^K(s-)=x, \theta \leq \bfk_y \tau_\bfk(x, \zeta^K(s-))}(\bfe(y) - \bfe(x)) Q_{x, \bfk, y}(ds, d\theta), \\
		\zeta^K(t) &= \overline{\bfz}^K + \frac{1}{K} \sum_{(x, \bfk, y) \in J^*} h_1(x, \bfk, y) \int_0^t \int_0^{+\infty} \setind{\theta \leq K \zeta^K_{1,x}(s-) \bfk_y \tau_\bfk(x, \zeta^K(s-))} Q_{x, \bfk, y}(ds, d\theta) \\
		&+ \frac{1}{K} \sum_{(x, \bfk) \in J}  h_0(x, \bfk) \int_0^t \int_0^{+\infty} \setind{\theta \leq K \zeta^K_{0,x}(s-)  \tau_{\bfk}(x,\zeta^K(s-))} Q_{x, \bfk}(ds, d\theta).
	\end{aligned}
	\end{equation}
}

With these notations, we establish the following control on the first moment of the population size.
\begin{lemma}
\label{lemma:moments}

Let $t > 0$. Under Assumption \ref{hyp:lipschitz}, there exists $C(t, \bfz) > 0$ such that for any $x \in \X$,
\[ \sup_{K \geq 1} \E_{(x, \overline{\bfz}^K)}\Big[\sup_{s \in [0,t]} \| \zeta^K \|_1\Big] \leq C(t, \bfz) < \infty. \]
\end{lemma}

\begin{proof}
	Let $n_0 = \sup_K \| \overline{\bfz}_K \|_1$, which is finite by definition of $\bfz_K$. Fix $K \geq 1$. The key lies in noticing that the number of particles alive at time $t$ in $\zeta^K$ can be dominated almost surely by a single type branching process $V$ coupled to $\zeta^K$ starting from $V(0) = n_0$, where each particle dies and is replaced by $n > 0$ particles at rate 
	\[ b_n = n \sum_{\bfk : \| \bfk \|_1 = n} \| \tau_\bfk \|_\infty.\]
In particular, as $J$ is a finite set, it holds that $\sum_n b_n < \infty$. Absence of deaths further implies that $\E[\sup_{s \leq t} V(s)] = \E[V(t)] < \infty$.
Thus for any $K \geq 1$, 
\[ \E_{\overline{\bfz}^K}\Big[ \sup_{s \in [0,t]}\| \zeta^K(s)\|_1 \Big] \leq \E\Big[V(t)\Big], \]
which concludes. 
\end{proof}

\paragraph{Quantification of the speed of convergence.} Let us start by establishing the following Lemma which quantifies the approximation error of the spinal population process by its large population limit. 

\begin{lemma} 
\label{lem:lln-error}
Let $t > 0$. Under Assumption \ref{hyp:lipschitz}, there exists $C(t, \bfz) > 0$ such that for any $x \in \X$, for any $K \geq 1$ and $\varepsilon_K > 0$, 
	\begin{equation*}
 	\P_{(x, \overline{\bfz}^K)}\left(\sup_{s \in [0,t]} \|\zeta^K(s) - z(s) \|_1\geq \varepsilon_K \right) \leq \frac{C(t, \bfz)}{\sqrt{K} \varepsilon_K}. 
	\end{equation*}
\end{lemma}

\begin{proof}
Let us define the large population limit $z$ on the same state space $\X^*$ as $\zeta^K$. For any $t \geq 0$,
\begin{equation}
\label{eq:lln-sde}
	z(t) = \overline{\bfz} + \sum_{(x,\bfk) \in J} h_0(x, \bfk) \int_0^t z_{0,x}(s) \tau_\bfk(x, z(s)) ds,
\end{equation}
with $\overline{\bfz} = \sum_{x \in \X} \bfz_x \bfe(0,x)$. In particular, $\bfz_{1,x}(t) = 0$ for any $x \in \X$ and $t \geq 0$, as expected since the spine becomes negligible as population size grows large. 

In order to control our quantity of interest $\|\zeta^K - z\|_1$, we introduce some notation. For $(x, \bfk) \in J$, let 
\begin{equation*}
\widetilde{Q}_{x,\bfk}(ds, d\theta) = Q_{x,\bfk}(ds, d\theta) - ds \, d\theta
\end{equation*}
be the compensated martingale-measure associated to $Q_{x, \bfk}$, and consider the following martingale:
\begin{equation*}
M^K_{x,\bfk}(t) = \frac{\|h_0(x, \bfk)\|_1}{K} \int_0^t \int_0^{+\infty} \setind{\theta \leq K \zeta^K_{0,x}(s-)\tau_\bfk(x, \zeta^K(s-))}\widetilde{Q}_{x,\bfk}(ds, d\theta).
\end{equation*}

Finally, let 
\begin{equation*}
\begin{aligned} 
M^K(t) &= \sum_{(x, \bfk) \in J} M^K_{x, \bfk}(t), \\
A^K(t) &= \Big\| \overline{\bfz}^K_0 - \overline{\bfz}_0 + \sum_{(x, \bfk) \in J} h_0(x, \bfk) \int_0^t (\zeta^K_{0,x}(s) \tau_\bfk(x, \zeta^K(s)) - z_{0,x}(s) \tau_\bfk(x, z(s)) ) ds \Big\|_1, \\
B^K(t) &= \Big\| \frac{1}{K} \sum_{(x, \bfk, y) \in J^*} h_1(x, \bfk, y) \int_0^t \int_0^{+\infty} \setind{\theta \leq K \zeta^K_{1,x}(s-) \bfk_y \tau_\bfk(x, \zeta^K(s-))} Q_{x, \bfk, y}(ds, d\theta) \Big\|_1.
\end{aligned}
\end{equation*}

With these notations, Equations \eqref{eq:psi-spine} and \eqref{eq:lln-sde} ensure that, for any $t \geq 0$ 
\begin{equation*}
\|\zeta^K(t) - z(t)\|_1 \leq M^K(t) + A^K(t) + B^K(t).
\end{equation*}
We thus turn to controlling $M^K$, $A^K$ and $B^K$ over any given interval $[0,t]$.  

Start by noticing that, for any $(x, \bfk) \in J$, $M^K_{x, \bfk}$ is square integrable. Indeed, as $\tau_{\bfk}$ is bounded and using Lemma \ref{lemma:moments}, there exists a positive constant $c_{x, \bfk}(t, \bfz)$ such that, for any $K \geq 1$,
\begin{equation*}
\E\left[\int_0^t  \int_0^{+\infty} \left( \frac{\|h_0(x, \bfk)\|_1}{K} \setind{\theta \leq K \zeta^K_{0,x}(s)\tau_\bfk(x, \zeta^K(s))} \right)^2 d\theta ds \right] \leq \frac{c_{x, \bfk}(t, \bfz)}{K}.
\end{equation*}
It follows that its quadratic variation is given by 
\begin{equation*}
\langle M^K_{x, \bfk} \rangle (t) =  \int_0^t  \frac{\|h_0(x, \bfk)\|^2_1}{K} \zeta^K_{0,x}(s)\tau_\bfk(x, \zeta^K(s)) ds.
\end{equation*}

As the family of Poisson Point Processes $(Q_{x, \bfk}, (x, \bfk) \in J)$ is independent, $M^K$ thus is itself a square integrable martingale of quadratic variation equal to 
\[ \langle M^K \rangle(t) = \sum_{(x, \bfk) \in J} \langle M^K_{x, \bfk} \rangle(t),\] 
and Doob's inequality further leads to 
\begin{equation}
\label{eq:lln-aux1}
\E \left[ \sup_{s \in [0,t]} (M^K(s))^2 \right] \leq 4 \sum_{(x, \bfk) \in J} \E[\langle M^K_{x, \bfk} \rangle(t)].
\end{equation}

Since $J$ is a finite set, it follows that there exists a positive, finite constant $c(t, \bfz)$ such that, for any $K \geq 1$,
\begin{equation*}
	\E \left[\sup_{s \in [0,t]} (M^K(s))^2 \right] \leq \frac{c(t, \bfz)}{K}.
\end{equation*}
Noticing that $M^K \geq 0$ almost surely, Cauchy-Schwarz inequality finally leads to 
\begin{equation*}
	\E \left[\sup_{s \in [0,t]} M^K(s) \right] \leq \sqrt{\frac{c(t, \bfz)}{K}}.
\end{equation*}

Let us now turn to $A^K$. By definition, there exists a finite constant $c(\bfz) > 0$ such that, for any $K \geq 1$,
\begin{equation*}
\| \overline{\bfz}^K - \overline{\bfz} \|_1 \leq \frac{c(\bfz)}{K}. 
\end{equation*}
Recall that the set $J$ is finite and $\tau_\bfk(x, \cdot)$ is Lipschitz-continuous according to Assumption \ref{hyp:lipschitz}. In addition, by continuity, $s \mapsto \| z(s) \|_1$ is bounded on $[0,t]$. As a consequence, there exists $c' (t, \bfz) > 0$ such that, for any $K \geq 1$, 
\begin{equation*}
\begin{aligned}
	\sup_{s \in [0,t]} A^K(s) & \leq \frac{c(\bfz)}{K} + \sum_{(x,\bfk) \in J} \| h_0(x, \bfk) \|_1 \int_0^t | \zeta^K_{0,x}(s) \tau_\bfk(x, \zeta^K(s)) - z_{0,x}(s) \tau_\bfk(x, z(s)) | ds \\
	& \leq \frac{c(\bfz)}{K} + c'(t, \bfz)  \int_0^t \| \zeta^K(s) - z(s) \|_1 ds.
\end{aligned}
\end{equation*}
In particular, we obtain that for any $K \geq 1$,
\begin{equation}
\label{eq:lln-aux2}
	\E \left[ \sup_{s \in [0,t]} A^K(s) \right] \leq \frac{c(\bfz)}{K} + c'(t, \bfz) \int_0^t \E \left[ \sup_{u \in [0,s]} \|\zeta^K(u) - z(u) \| ds \right].
\end{equation}

Finally, notice that it follows from Equation \eqref{eq:psi-spine} that for any $t \geq 0$ and $x \in \X$, we have $\zeta^K_{1,x} \leq K^{-1}$ almost surely. As further the set $J^*$ is finite and the reproduction rates bounded, there exists a finite, non-negative constant $c$ such that, for any $K \geq 1$, 
\begin{equation}
\label{eq:lln-aux3}
\E\left[ \sup_{s \in [0,t]} B^K(s) \right] \leq \sum_{(x,\bfk,y) \in J^*} \|h_1(x, \bfk, y) \|_1 \E\left[ \int_0^t \zeta^K_{1,x}(s) \bfk_y \tau_\bfk(x, \zeta^K(s)) ds \right] \leq \frac{tc}{K}.
\end{equation}

Combining Equations \eqref{eq:lln-aux1}, \eqref{eq:lln-aux2} and \eqref{eq:lln-aux3} yields the existence of finite, non-negative constants $c_1(t, \bfz)$ and $c_2(t, \bfz)$ such that, for any $K \geq 1$,
\begin{equation*}
	\E \left[ \sup_{s \in [0,t]} \|\zeta^K(s) - z(s) \|_1 \right] \leq \frac{c_1(t, \bfz)}{\sqrt{K}} + c_2(t, \bfz) \int_0^t \E \left[ \sup_{u \in [0,s]} \|\zeta^K(u) - z(u) \|_1 \right] ds.
\end{equation*}
Hence, by Gronwall's Lemma, there exists $C(t, \bfz)$ such that, for any $K \geq 1$, 
\begin{equation*}
	\E \left[ \sup_{s \in [0,t]} \|\zeta^K(s) - z(s) \|_1 \right] \leq \frac{C(t, \bfz)}{\sqrt{K}}.
\end{equation*}
The conclusion follows from Markov's inequality.
\end{proof}

Next, we focus on quantifying how well $\Upsilon$ approaches $Y^K$, yielding the approximation error associated to replacing $(Y^K, \zeta^K)$ by $(\Upsilon, z)$ on finite time intervals. In order to achieve this, we couple $\Upsilon$ and $Y^K$, by defining $\Upsilon$ as the unique strong solution to an SDE driven by the same family of Poisson Point Processes $(Q_{x, \bfk, y}, (x, \bfk, y) \in J^*)$ as in Equation \eqref{eq:psi-spine}:
\begin{equation*}
\Upsilon(t) = \bfe(y_0) + \sum_{(x, \bfk, y) \in J^*} \int_0^t \int_0^{+\infty} \setind{\Upsilon(s-)=x, \theta \leq \bfk_y \tau_\bfk(x, z(s))}(\bfe(y) - \bfe(x)) Q_{x, \bfk, y}(ds, d\theta).
\end{equation*}
We are now ready to state our result.

\begin{lemma}
\label{thm:error}
	Let $t > 0$. Under Assumption \ref{hyp:lipschitz}, there exists a constant $C(t, \bfz) > 0$
	such that for every $K \geq 1$ and $\varepsilon_K > 0$, letting $\alpha_K = (\sqrt{K} \varepsilon_K)^{-1}$,
	\begin{equation*}
	\P(\forall s \in [0,t], Y^K(s) = \Upsilon(s) \text{ and } \|\zeta^K(s) - z(s)\|_1 \leq \varepsilon_K) \geq 1 - C(t, \bfz)(\alpha_K + \varepsilon_K).
	\end{equation*}
\end{lemma}

\begin{proof}

Consider Lemma \ref{lemma:error} with $(Y_1, \zeta_1) = (Y^K, \zeta^K)$ and $(Y_2, \zeta_2) = (\Upsilon, z)$. Thanks to Lemma \ref{lem:lln-error}, we thus obtain that 
	\begin{equation}
	\label{eq:precision-spine}
		\P(\exists t \leq t : Y^K(t) \neq \Upsilon(t)) \leq c(t, \bfz) \alpha_K + C t \varepsilon_K.
	\end{equation}
The conclusion follows by combining Lemma \ref{lem:lln-error} and Equation \eqref{eq:precision-spine}.

\end{proof}

In order to simplify notations, for any function $G : \mathbb{D}([0,t], \X)  \times \mathbb{D}([0,t], \R_+^D) \to \R$ as well as càdlàg trajectories $(x(s))_{s \leq t}$ in $\X$ and $(\bfz(s))_{s \leq t} = (\bfz_{i,x}(s))_{(i,x) \in \X^*, s \leq t}$ in $\R_+^{2D}$, we let 
\begin{equation*}
	G((x(s), \bfz(s))_{s \leq t}) = G((x(s), \proj{\bfz}(s))_{s \leq t}).
\end{equation*}
We are now ready to establish the main result. 

\begin{proof}[Proof of Proposition \ref{prop:lln-sample}] 
Throughout the proof, we write $\bbY^K = (Y^K, \zeta^K)$ for the spinal construction in the initial population process, and $\bbY = (\Upsilon, z)$ for the spinal process in the large population limit. 

Let $t \geq 0$ and $F \in \mathfrak{L}$. It follows from Theorem 1 in \cite{bansayeSpineInteractingPopulations2024} that 
\begin{equation*}
	\E_{\bfz^K}\left[\sum_{\substack{u \in \generation{t},\\ u \succeq u_x}} F((x_u(s), Z^K(s))_{s \leq t})\right] = \E \left[H((\bbY^K(s))_{s \leq t}) \right],
\end{equation*}
where 
\begin{equation*}
H((\bbY^K(s))_{s \leq t}) = \mathrm{e}^{\int_0^t \lambda(\bbY^K(s))ds}F((\bbY^K(s))_{s \leq t}).
\end{equation*}

Consider the event 
\begin{equation*}
	A_K = \left\{ \forall s \in [0,t], Y^K(s) = \Upsilon(s) \text{ and } \|\zeta^K(s) - z(s)\|_1 \leq \varepsilon_K \right\}.
\end{equation*}
We have 
\begin{equation*}
\begin{aligned}
	\E[|H((\bbY^K(s))_{s \leq t}) - H((\bbY(s))_{s \leq t}) |] & \leq \E[|H((\bbY^K(s))_{s \leq t}) - H((\bbY(s))_{s \leq t}) | \ind_{A_K}] \\ 
	& + \E[|H((\bbY^K(s))_{s \leq t}) - H((\bbY(s))_{s \leq t}) | \ind_{A_K^C}]
\end{aligned}
\end{equation*}

On the one hand, it follows from the definition of $A_K$ and continuity of $s \mapsto z(s)$ that there exists a compact subset $\Z \in \R_+^D$ such that for every $K$,  
\begin{equation}
\label{eq:final-bound}
	 \P\big(\forall s \in [0,t], \; \zeta_K(s) \in \Z \big| A_K \big) = 1. 
\end{equation}
As $J$ is a finite set and $F \in \mathfrak{L}$, Assumption \ref{hyp:lipschitz} implies that there exists a constant $M$ depending on $t$ and $F$ such that for any $x \in \mathbb{D}([0,t],\X)$ and $K \geq 1$, on the event $A_K$,
\begin{equation*}
\label{eq:temp2}
	|H((x(s), \zeta^K(s))_{s \leq t}) - H((x(s), z(s))_{s \leq t}) | \leq M \sup_{s \in [0,t]} \|\zeta^K(s) - z(s)\|_1.
\end{equation*}
It follows that
\begin{equation*}
\E[|H((\bbY^K(s))_{s \leq t}) - H(\bbY(s))_{s \leq t}) | \ind_{A_K}] \leq M \varepsilon_K.
\end{equation*}
On the other hand, the boundedness of $H$ and Lemma \ref{thm:error} imply the existence of a constant $c$ such that 
\begin{equation*}
\E[|H((\bbY^K(s))_{s \leq t}) - H((\bbY(s))_{s \leq t}) | \ind_{A_K^C}] \leq c(\alpha_K + \varepsilon_K).
\end{equation*}
Taken together, we thus obtain existence of a constant $C$ such that 
\begin{equation*}
\E[|H((\bbY^K(s))_{s \leq t}) - H((\bbY(s))_{s \leq t})|] \leq C (\alpha_K + \varepsilon_K).
\end{equation*}
This concludes the proof. 
\end{proof}

\section{Discussion}
\label{sec:discussion}

First, we have introduced a time-inhomogeneous spinal process allowing to gain insight on the survivorship bias associated to any sampling weight $\psi$. 

Indeed, the corresponding many-to-one formula does not require stochastic weighting of trajectories. This implies that the bias of reproduction rates relying on the function $m$ defined in Equation \eqref{eq:mpsi} accurately depicts the survivorship bias. In addition, the stochastic exponential weight associated to the many-to-one formula for the time-homogeneous spinal process implies that rare trajectories may have tremendous impact, making Monte-Carlo estimations delicate. As a consequence, the time-inhomogeneous spinal process may facilitate the numerical evaluation of the many-to-one formula. However, due to the time-inhomogeneity, 
simulating trajectories of the time-inhomogeneous spinal process through standard algorithms may be expensive in terms of computation time \cite{thanhSimulationBiochemicalReactions2015}.

A desirable extension of Theorem \ref{thm:spine} consists in capturing the whole tree, instead of being restricted to the type evolution along sampled lineages. We expect this to be achievable, using a classical induction argument \cite[Theorem 1]{bansayeSpineInteractingPopulations2024}. In addition, one may be interested in sampling more than one individual, in the spirit of many-to-one formulas for forks in branching processes \cite{marguetUniformSamplingStructured2019}. A possible strategy for achieving this would be a double spine construction. More precisely, one may augment the type space of the spinal process to distinguish spinal and non-spinal individuals, and then consider the spinal process associated to (re-)sampling in this population process. 

Second, under appropriate assumptions, we have focused on sampling in the deterministic large population limit, and quantified the associated approximation error. In particular, the process describing the spine's offspring then corresponds to a time-inhomogeneous multi-type branching process whose reproduction rates depend on a changing environment, given by $z$. If we assume that $z$ admits a stable equilibrium, then starting from this equilibrium, the offspring is described by a classical multi-type branching process \cite{calvezDynamicsLineagesAdaptation2022}. 

A natural perspective of our work is to extend our results on the empirical distribution of ancestral lineages in the large population limit over longer time scales. In order to achieve this, we need some control of the form:  
\begin{equation*}
	\P\left(\sup_{s \in [0, t_K]} \|\zeta^K(s) - z(s)\|_1 \geq \varepsilon_K\right) \leq \alpha_K,
\end{equation*}
with $\varepsilon_K$ and $\alpha_K$ converging to zero and $t_K$ growing to infinity, as $K$ grows large. We consider expecting such control to be reasonable. Indeed, it corresponds to understanding and controlling the fluctuations of the finite-population process $\zeta^K$ around its deterministic limit, which is a well studied question with several classical regimes: Gaussian fluctuations for $\varepsilon_K = O(K^{-1/2})$, which are related to the diffusion approximation, moderate deviations with $\varepsilon_K = O(K^{-p})$ for $p \in (0,1/2)$ and large deviations where $\varepsilon_K = O(1)$ \cite{brittonStochasticEpidemicsHomogeneous2019, pardouxModerateDeviationsExtinction2020, prodhommeStrongGaussianApproximation2023}. In particular, moderate and large deviations appear to be interesting regimes, as they allow to consider longer time scales. Nevertheless, in the case of the spinal population process, they are not immediate, due to a boundary problem arising from the fact that the spinal individual becomes negligible in large populations. As a consequence, these considerations are left for a futur work.

\paragraph{Acknowledgments.} The author thanks Vincent Bansaye for stimulating discussions and his continuous support during the elaboration of this work. The author also is grateful to Aline Marguet, Charline Smadi and Charles Medous for fruitful discussions. The author thanks an anonymous reviewer for very constructive feedback. This work was partially funded by the Chair “Modélisation Mathématique et Biodiversité” of VEOLIA- Ecole Polytechnique-MNHN-F.X.

\bibliographystyle{plain} 
\bibliography{spine}

\begin{appendix}

\section{Proof of Proposition \ref{prop:psi-process}}
\label{appdx:existence}

The proof is decomposed in three steps, establishing (i) existence and (ii) uniqueness of the solution to Equation \eqref{eq:def-sde} by classical arguments, before (iii) characterizing the associated semi-group $R^{(t)}$. For ease of notation, throughout the proof, for $0 \leq s \leq t$ we let $\couplespine(s) = (Y^{(t)}(s), \zeta^{(t)}(s))$.
\medskip

(i) \emph{Existence.} First, notice that by Assumptions \ref{hyp:avg} and \ref{hyp:bound-m}, both $\rho^{(t)}_{y, \bfk}$ and $\hat{\rho}^{(t)}_\bfk$ are bounded for any $y \in \X$ and $\bfk \in \Z$. As a consequence, existence of at least one solution to Equation \eqref{eq:def-sde} is ensured, as the associated sequence of jump times $(T_k)_{k \geq 0}$ cannot admit an accumulation point on $[0,t]$. 
\medskip

(ii) \emph{Uniqueness.} Subsequently, in order to establish uniqueness, let us show by induction that for any $k \geq 0$ such that $T_k \leq t$, $(T_k, \couplespine(T_k))$ is entirely determined by $(\couplespine(0), Q, \widehat{Q})$. As $T_0 = 0$, initialization of the induction argument is immediate. If the property holds for $k \geq 1$, then by construction, $T_{k+1}$ only depends on $(T_k, \couplespine(T_k), Q, \widehat{Q})$. Similarly, given $T_{k+1}$ and the corresponding atoms $A_{k+1}$ and $\widehat{A}_{k+1}$ of $Q$ and $\widehat{Q}$, it is clear that $\couplespine(T_{k+1})$ is fixed by $(T_{k+1}, A_{k+1}, \widehat{A}_{k+1}, \couplespine(T_k))$. The desired conclusion thus is a consequence of the induction hypothesis. 
\medskip

(ii) \emph{Characterization of $R^{(t)}$}. In order to establish Equation \eqref{eq:defR}, it is sufficient to show that for any bounded function $f$ on $\configS$ and $(x, \bfz) \in \configS$, the function 
\begin{equation*}
	\tau \mapsto R^{(t)}_{s, \tau}f(x,\bfz) = \E[f(\couplespine(\tau)) | \couplespine(s) = (x, \bfz)]
\end{equation*}
is right differentiable at $\tau = s$. Indeed, it then follows that Equation \eqref{eq:defR} holds with the operator $A^{(t)}$ defined by
\begin{equation}
\label{eq:diffR}
	\forall f : \configS_K \to \R_+\; \forall (x,\bfz) \in \configS_K, \quad A^{(t)}_s f(x, \bfz) = \lim_{h \to 0+} \frac{1}{h} \left( R^{(t)}_{s,s+h}f(x, \bfz) - f(x, \bfz) \right).
\end{equation}
As we will see, computing the right-hand side of Equation \eqref{eq:diffR} leads to Equation \eqref{eq:defA}.

Let $f : \configS \to \R_+$ and $(x, \bfz) \in \configS$. We introduce the following notations. For any $(y, \bfk) \in \configS$ such that $\tau_\bfk(x, \bfz) > 0$,
\begin{equation*}
 	\mathfrak{d}_{y, \bfk} f(x, \bfz) = f(y, \bfz + \bfk - \bfe(x)) - f(x, \bfz).
\end{equation*}
Further, for any $y \in \X$ such that $\bfz_y > 0$ and $\bfk \in \Z$ such that $\tau_\bfk(y, \bfz) > 0$, let 
\begin{equation*}
	\widehat{\mathfrak{d}}_{y, \bfk} f(x, \bfz) = f(x, \bfz + \bfk - \bfe(y)) - f(x, \bfz).
\end{equation*}

Equation \eqref{eq:def-sde} then ensures that, on the event $\couplespine(s) = (x,\bfz)$, we have for any $h \in [0,t-s]$: 
\begin{equation*}
\begin{aligned}
	f(\couplespine(&s+h)) - f(x, \bfz) = \int_s^{s+h} \int_{E} \setind{\theta \leq \rho^{(t)}_{y, \bfk}(r, \couplespine(r-))} \mathfrak{d}_{y, \bfk} f(\couplespine(r-)) Q(dr, d\theta, n(dy, d\bfk)) \\
	&+ \int_s^{s+h} \int_{E} \setind{\theta \leq (\zeta^{(t)}_y(r-) - \setind{Y^{(t)}(r-) = y}) \widehat{\rho}^{(t)}_{\bfk}(r, y,\couplespine(r-))} \widehat{\mathfrak{d}}_{y, \bfk} f(\couplespine(r-)) \widehat{Q}(dr, d\theta, n(dy, d\bfk)).
\end{aligned}
\end{equation*}

Notice that, for instance,
\begin{equation*}
\begin{aligned}
\E[\int_s^{s+h} \int_{E} &\setind{\theta \leq \rho^{(t)}_{y, \bfk}(r, \couplespine(r-))}\mathfrak{d}_{y, \bfk} f(\couplespine(r-)) Q(dr, d\theta, n(dy, d\bfk)) | \couplespine(s)=(x, \bfz)] \\
&= \E[\int_s^{s+h} \sum_{(y, \bfk) \in \configS} \rho^{(t)}_{y, \bfk}(r, \couplespine(r)) \mathfrak{d}_{y, \bfk} f(\couplespine(r)) dr | \couplespine(s)=(x, \bfz)].
\end{aligned}
\end{equation*}

On the one hand, almost surely, 
\begin{equation*}
\begin{aligned}
\lim_{h \to 0+} \frac{1}{h} \int_s^{s+h} \sum_{(y, \bfk) \in \configS} & \rho^{(t)}_{y, \bfk}(r,\couplespine(r)) \mathfrak{d}_{y, \bfk} f(\couplespine(r)) dr = \sum_{(y, \bfk) \in \configS} \rho^{(t)}_{y, \bfk}(s,\couplespine(s)) \mathfrak{d}_{y, \bfk} f(\couplespine(s)).
\end{aligned}
\end{equation*}

By Assumptions \ref{hyp:avg} and \ref{hyp:bound-m}, $\sum_{y, \bfk} \|\rho^{(t)}_{y, \bfk} \|_\infty < \infty.$
Thus, there exists a constant $C > 0$ such that for any $h \in [0, t-s]$,
\begin{equation*}
\frac{1}{h} \int_s^{s+h} \left| \sum_{(y, \bfk) \in \configS} \rho^{(t)}_{y, \bfk}(r, \couplespine(r)) \mathfrak{d}_{y, \bfk} f(\couplespine(r)) \right| dr \leq C \| f \|_\infty < \infty.
\end{equation*}
Taken together, we obtain by dominated convergence:
\begin{equation*}
\begin{aligned}
\frac{1}{h} \E[\int_s^{s+h} \!\!\! \int_{E} \!\! \setind{\theta \leq \rho^{(t)}_{y, \bfk}(r, \couplespine(r-))} \mathfrak{d}_{y, \bfk} f(\couplespine(r-)) & Q(dr, d\theta, n(dy, d\bfk)) | \couplespine(s)=(x, \bfz)] \\
& \xrightarrow[h \to 0+]{} \sum_{(y, \bfk) \in \configS} \rho^{(t)}_{y, \bfk}(s, x, \bfz)  \mathfrak{d}_{y, \bfk} f(x, \bfz).
\end{aligned}
\end{equation*}

The other terms arising on the right-hand side of Equation \eqref{eq:diffR} can be treated analogously. This leads to the desired result.

\section{Proof of Lemma \ref{lem:dsmpsi}} 
\label{appdx:m}

Let $(x, \bfz) \in \configS$, $t \geq 0$ and $h > 0$. The Markov property ensures that 
\begin{equation*}
	m(x, \bfz, t+h) = \popE[\sum_{\substack{u \in \generation{h}\\u \succeq u_x(0)}} m(x_u(h), Z(h), t)].
\end{equation*}
For $i \geq 1$, Let $T_i$ be the time of the $i$-th jump of the population process. Then on the one hand, if $T_1 > h$, then the population at time $h$ is identical to the population at time $0$, and thus:
\begin{equation*}
	m(x, \bfz, t+h) = \popE[\sum_{\substack{u \in \generation{h}\\u \succeq u_x(0)}} m(x_u(h), Z(h), t) \setind{T_1 < h}] + m(x, \bfz, t) \P_\bfz(T_1 > h).
\end{equation*}
Similarly, on the event $\{T_1 < h < T_2\}$, $Z(h) = Z(T_1)$ whence
\begin{equation*}
	\popE[\sum_{\substack{u \in \generation{h}\\u \succeq u_x(0)}} m(x_u(h), Z(h), t) \setind{T_1 < h}] = a(h) + b(h),
\end{equation*}
where
\begin{equation*}
\begin{aligned}
	a(h) &= \popE[\sum_{\substack{u \in \generation{h}\\u \succeq u_x(0)}} m(x_u(T_1), Z(T_1), t) \setind{T_1 < h < T_2}], \\
	b(h) &= \popE[\sum_{\substack{u \in \generation{h}\\u \succeq u_x(0)}} m(x_u(h), Z(h), t) \setind{T_2 < h}].
\end{aligned}
\end{equation*}
As a consequence, we obtain that 
\begin{equation*}
	m(x, \bfz, t+h) - m(x, \bfz, t) = A(h) + B(h),
\end{equation*}
with $A(h) = a(h) - m(x, \bfz, t)\P_\bfz(T_1 < h < T_2)$ and $B(h) = b(h) - m(x, \bfz, t)\P_\bfz(T_2 < h)$.
\medskip

Let us first focus on $B(h)$. For $\bfv \in \Z_K$, let us write $\Lambda(\bfv) = \sum_{y \in \X} \sum_{\bfk \in \Z_K} \bfv_y \tau_\bfk(y, \bfv)$ for the total jump rate in a population whose type distribution is given by $\bfv$. In particular, $\Lambda$ is bounded on $\Z_K$. Using the law of $T_1$ given $Z(0) = \bfz$ and the law of $T_2 - T_1$ given $T_1$ and $Z(T_1)$, we then obtain:
\begin{equation*}
	\P_\bfz(T_2 < h) = \int_0^h e^{-\Lambda(\bfz)t_1} \sum_{\substack{y \in \X\\\bfk \in \Z}} \bfz_y \tau_\bfk(y, \bfz) \int_0^{h-t_1}\Lambda(\bfz + \bfk - \bfe(y))\e^{-\Lambda(\bfz + \bfk - \bfe(y))t_2} dt_2 dt_1.
\end{equation*}
Assumption \ref{hyp:avg} thus leads to the existence of a positive constant $C(\bfz)$ such that 
\begin{equation}
\label{eq:T2}
	\P_\bfz(T_2 < h) \leq C(\bfz) h^2.
\end{equation}
Thus $m(x, \bfz, t)\P_\bfz(T_2 < h)$ is of order $o(h)$, and it remains to focus on $b(h)$. According to the Markov property,
\begin{equation*}
	b(h) = \E_\bfz[ \sum_{\substack{u \in \generation(t+h)\\u \succeq u_x(0)}}  \setind{T_2 < h}].
\end{equation*}

The key lies in noticing that the number of children alive at time $t + h$ and who descend from $u_x(0)$ can be dominated almost surely by a single type branching process $V$ coupled to $Z$ starting from $V(0) = 1$, where each particle dies and is replaced by $n > 0$ particles at rate 
\[ \sum_{k \in \Z : \| k\|_1 = n} \| \tau_\bfk \|_\infty .\] 
Assume that $\sum_{\bfk} \| \bfk \|_1^\alpha \| \tau_\bfk \|_\infty < \infty$ holds for $\alpha > 2$, which implies that $\E[V^{\alpha}] < \infty$. Holder's inequality and Equation \eqref{eq:T2} lead to the existence of a positive constant $C'(\bfz)$ such that
\[b(h) \leq  \E[V^{\alpha}]^{1/\alpha} C'(\bfz) h^{2(\alpha-1)/\alpha}.\] 
In particular, since $\alpha > 2$, it follows that $2(\alpha-1)/\alpha > 1$. We thus conclude that $b(h) = o(h)$. 
\medskip

Let us now focus on $A(h)$. Using the law of $T_1$ given $Z(0) = \bfz$, we have 
{\small
\begin{equation*}
\begin{aligned}
	A(h) &= \int_0^h e^{-\Lambda(\bfz)t_1} \e^{-\Lambda(\bfz + \bfk - \bfe(x))(h-t_1)} dt_1 \sum_{\bfk \in \Z} \tau_\bfk(x, \bfz) \sum_{y \in \X} \bfk_y (m(y, \bfz + \bfk - \bfe(x), t) - m(x, \bfz,t)) \\
	& + \sum_{\substack{y \in \X \\ \bfk \in \Z}} \int_0^h e^{-\Lambda(\bfz)t_1} \e^{-\Lambda(\bfz + \bfk - \bfe(y))(h-t_1)} dt_1 \, \tau_\bfk(y,\bfz) (m(x, \bfz + \bfk - \bfe(y), t) - m(x, \bfz,t)),
\end{aligned}
\end{equation*}
}
from which it follows that
\begin{equation}
\label{eq:Ah}
\begin{aligned}
	\frac{A(h)}{h} \xrightarrow[h \to 0^{+}]{} & \sum_{\bfk \in \Z} \tau_\bfk(x, \bfz) \left(\sum_{y \in \X} \bfk_y m(y, \bfz + \bfk - \bfe(x), t) - m(x, \bfz,t) \right) \\
	& + \sum_{\substack{y \in \X \\ \bfk \in \Z}} (\bfz_y - \setind{x = y})\tau_\bfk(y,\bfz) (m(x, \bfz + \bfk - \bfe(y), t) - m(x, \bfz,t)).
\end{aligned}
\end{equation}
As a consequence, right differentiability of $t \mapsto m(x, \bfz, t)$ is established, and its right derivative is given by the right-hand side of Equation \eqref{eq:Ah}. As this corresponds to a continuous function on $\R_+$, we deduce that $t \mapsto m(x, \bfz, t)$ is differentiable on $\R_+$ (see \emph{e.g.} Corollary 1.2 of Chapter 2 in \cite{pazySemigroupsLinearOperators2012}) and 
\begin{equation*}
\begin{aligned}
	\frac{d}{dt} m(x, \bfz, t) & = \sum_{\bfk \in \Z} \tau_\bfk(x, \bfz) \left(\sum_{y \in \X} \bfk_y m(y, \bfz + \bfk - \bfe(x), t) - m(x, \bfz,t) \right) \\
	& + \sum_{\substack{y \in \X \\ \bfk \in \Z}} (\bfz_y - \setind{x = y}) \tau_\bfk(y,\bfz) (m(x, \bfz + \bfk - \bfe(y), t) - m(x, \bfz,t)).
\end{aligned}
\end{equation*}
This concludes the proof.

\section{Detail on the proof of Proposition \ref{prop:spine-bounded}}
\label{appdx:monotone-class}

Throughout this proof, for readability, we will make use of the following notations. On the one hand, for $t \geq 0$ and $u \in \generation{t}$, let
\begin{equation*}
 	\couple_u(t) = (x_u(t), Z(t)).
\end{equation*}
Similarly, for $0 \leq s \leq t$, we let
\begin{equation*}
	\couplespine(s) = (Y^{(t)}(s), \zeta^{(t)}(s)).
\end{equation*}

Let us start by showing that Equation \eqref{eq:spine} holds for $F((x(s), \bfz(s))_{s \leq t}) = \prod_{j=1}^k f_j(x(s_j), \bfz(s_j))$ where $k \geq 1$, $0 \leq s_1 \leq \dots \leq s_k \leq t$ and $f_1, \dots, f_k: \configS_K \to \R_+$. 

This part of the proof proceeds by induction. For $k \geq 1$, let $H_k$ be the property that for any $0 \leq s_1 \leq \dots \leq s_k \leq t$ and $f_1, \dots, f_k: \configS_K \to \R_+$, 
\begin{equation*}
\begin{aligned}
	\popE [\!\!\!\!\sum_{u \in \generation{t},\, u \succeq u_x(0)}\!\!\!\! \psi(\couple_u(t)) \prod_{j=1}^k f_j(\couple_u(s_j))] = m(x,\bfz,t)\spineE[\prod_{j=1}^k f_j(\couplespine(s_j))].
\end{aligned}
\end{equation*}

Let us turn our attention to the initialization step. As $\configS_K$ is a finite set, a semi-group acting on non-negative functions on $\configS_K$ is uniquely characterized by its generator. Thus Lemma \ref{lem:gen} implies that the semi-groups $P^{ (t)}$ and $R^{(t)}$ are identical. Hence for any $s \in [0,t]$ and $f : \configS_K \to \R_+$, Equation \eqref{eq:defp} becomes 
\begin{equation*}
	 \popE[\sum_{u \in \generation{t},\, u \succeq u_x(0)}\psi(\couple_u(t)) f(\couple_u(s))] = m(x, \bfz, t) R^{(t)}_{0,s}f(x, \bfz).
\end{equation*}
This exactly corresponds to $H_1$ by definition of $R^{(t)}$.

Suppose now that $H_{k-1}$ is true for $k > 1$, and let us show that $H_{k}$ follows. Consider functions $f_1, \dots, f_k: \configS_K \to \R_+$ and $0 \leq s_1 \leq \dots \leq s_k \leq t$. Notice that
{ \small
\begin{equation*}
\begin{aligned}
	\popE[\!\!\!\!\sum_{\substack{u \in \generation{t}\\u \succeq u_x(0)}}\!\!\!\! &\psi(\couple_u(t)) \prod_{j=1}^k f_j(\couple_u(s_j))] = \\
	&\popE [\!\!\!\! \sum_{\substack{u \in \generation{s_{k-1}}\\u \succeq u_x(0)}} \prod_{j=1}^{k-1} f_j(\couple_u(s_j)) \E[\!\!\!\!\!\! \sum_{\substack{v \in \generation{t}\\v \succeq u_{x_u(s_{k-1})}}}\!\!\!\! \psi(\couple_u(t)) f_k(\couple_v(s_k))|X(s_{k-1})=\mathfrak{X}(\couple_u(s_{k-1}))]].
\end{aligned}	
\end{equation*}
}
As for $H_1$, equality of $P^{ (t)}$ and $R^{(t)}$ leads to:
{ \small
\begin{equation*}
\begin{aligned}
	\popE &[\!\!\!\! \sum_{\substack{u \in \generation{t}\\u \succeq u_x(0)}}\!\!\!\! \psi(\couple_u(t)) \prod_{j=1}^k f_j(\couple_u(s_j))] = \\
	 & \popE \Big[\!\!\!\! \sum_{\substack{u \in \generation{s_{k-1}}\\u \succeq u_x(0)}} m(\couple_u(s_{k-1}), t-s_{k-1})  \prod_{j=1}^{k-1} f_j(\couple_u(s_j))  \E[f_k(\couplespine(s_k))|\couplespine(s_{k-1}) = \couple_u(s_{k-1}) ]\Big].
\end{aligned}	
\end{equation*}
}
Equation \eqref{eq:ttau} allows to rewrite this as:
{ \small
\begin{equation*}
\begin{aligned}
	\popE [\!\!\!\! \sum_{\substack{u \in \generation{t}\\u \succeq u_x(0)}}&\!\!\!\! \psi(\couple_u(t)) \prod_{j=1}^k f_j(\couple_u(s_j))] =\\
	& \popE \Big[\!\!\!\! \sum_{\substack{u \in \generation{t}\\u \succeq u_x(0)}} \!\!\!\! \psi(\couple_u(t)) \prod_{j=1}^{k-1} f_j(\couple_u(s_j))  \E[f_k(\couplespine(s_k))|\couplespine(s_{k-1}) = \couple_u(s_{k-1}) ]\Big].
\end{aligned}	
\end{equation*}
}
Finally, $H_{k-1}$ yields:
\begin{equation*}
\begin{aligned}
	\popE [\!\!\!\! \sum_{\substack{u \in \generation{t}\\u \succeq u_x(0)}}\!\!\!\! & \psi(\couple_u(t)) \prod_{j=1}^k f_j(\couple_u(s_j))] \\
	& = m(x, \bfz,t) \spineE[\prod_{j=1}^{k-1} f_j(\couplespine(s_j)) \E[f_k(\couplespine(s_k)) | \couplespine(s_{k-1})]] \\
	& = m(x, \bfz,t) \spineE[\prod_{j=1}^k f_j(\couplespine(s_j))].
\end{aligned}	
\end{equation*}
This concludes the induction argument.

In order to obtain the desired result, we will reason using the monotone class theorem. Let us introduce the set 
\begin{equation*}
	I = \left\{\bigcap_{j=1}^k \{x \in \D([0,t], \configS_K): x(s_j) \in B_j\}, k \in \N, s_j \in [0,t], B_j \in \calP(\configS_K) \right\}
\end{equation*}
where $\calP(\configS_K)$ is the set of subsets of $\configS_K$. The set $I$ is a $\pi$-system, which induces the Borel $\sigma$-algebra $\mathcal{B}(\D([0,t], \configS_K))$ on the Skorokhod space $\D([0,t], \configS_K)$ (Theorem 12.5 in \cite{billingsleyConvergenceProbabilityMeasures1999}).
Further, define 
\begin{equation*}
	M = \{\mathbf{B} \in \mathcal{B}(\D([0,t], \configS_K)):\; \text{Equation \eqref{eq:spine} is satisfied for } F = \ind_{\mathbf{B}} \}.
\end{equation*}
$M$ is a monotone class which contains $I$ according to our induction argument. It thus follows from the monotone class theorem that $M =  \mathcal{B}(\D([0,t], \configS_K))$. In other words, for any $\mathbf{B} \in \mathcal{B}(\D([0,t], \configS_K))$, Equation \eqref{eq:spine} is satisfied for $F = \ind_\mathbf{B}$. As a consequence, Equation \eqref{eq:spine} holds for any positive measurable function $F: \D([0,t], \configS_K) \to \R_+$ as there exists an increasing sequence of simple functions converging pointwise to $F$, from which the result follows by monotone convergence.

\section{Proof of Proposition \ref{prop:lln-spine}}
\label{proof:prop-lln-spine}

Consider a family of independent Poisson Point Processes $(Q_{j}, j \in J)$ of intensity the Lebesgue measure on $\R_+^2$. The process $(\Upsilon^{(t)}, z)$ is characterized as follows: for any $s \in [0,t]$, for any initial condition $(x, \bfz) \in \overline{\configS}$,
\begin{equation*}
\begin{aligned}
z'(s) & = z(s)A(z(s)), \;\; z(0) = \bfz, \\
\Upsilon^{(t)}(s) & = \bfe(x) + \sum_{\bfk, y \in J} \int_0^t \int_0^{+\infty} \setind{\theta \leq \rho^{(t)}_{\bfk, y}(\Upsilon^{(t)}(u-), z(u), u)}(\bfe(y) - \bfe(\Upsilon^{(t)}(u-))) Q_{\bfk, y}(du, d\theta).
\end{aligned}
\end{equation*}

In order to see that the process $(\Upsilon^{(t)}, z)$ is well defined on $[0,t]$, start by noticing that existence and uniqueness of $z \in \C^1(\R_+)$ follows from the Cauchy-Lipschitz theorem and Assumption \ref{hyp:lipschitz}.
It remains to focus on existence and uniqueness of $\Upsilon^{(t)}$, given $z$.  This can be achieved through similar arguments as for steps (i) and (ii) in the proof of Proposition \ref{prop:psi-process}.  Indeed, as the reproduction rates are bounded, the function $s \mapsto \m(x, z(s), t-s)$ is bounded from below by $\exp(-\| \tau_0 \|_\infty t)$ on $[0,t]$, and the reproduction rates $\rho^{(t)}_{\bfk, y}$ are also bounded.
\bigskip

We now turn to characterizing the associated semi-group $\calR^{(t)}$. It is clear that for any $0 \leq r \leq s \leq t$, $\calR^{(t)}_{r,s}$ is a linear operator on $\mathrm{L}^\infty(\overline{\configS})$. Further, 
\begin{equation*}
	\forall f \in \mathrm{L}^\infty(\overline{\configS}), \quad \| \calR^{(t)}_{r,s} f \|_\infty \leq \|f\|_\infty.
\end{equation*}
Thus, $\calR^{(t)}_{r,s} : \mathrm{L}^\infty(\overline{\configS}) \to \mathrm{L}^\infty(\overline{\configS})$ and its operator norm equals one, as $\calR^{(t)}_{r,s}1 = 1$. 

Let us finally compute the generator of $\calR^{(t)}$: we want to show that for any $f \in \mathcal{C}^1(\overline{\configS})$, $s \leq t$ and $(x, \bfz) \in \overline{\configS}$,
\begin{equation*}
	\lim_{h \to 0+} \frac{1}{h} \left( \calR^{(t)}_{s,s+h}f(x,\bfz) - f(x,\bfz) \right) = \A^{(t)}_s f(x,\bfz).
\end{equation*}

Fix $s$ and define $T_1$ and $T_2$ as the times of the first and second reproduction events after time $s$. By definition, we have 
\begin{equation*}
	\calR^{(t)}_{s,s+h}f(x,\bfz) = A(h) + B(h) + C(h),
\end{equation*}
where 
\begin{equation*}
\begin{aligned}
	A(h) &= \E[f(\Upsilon^{(t)}(s+h), z(s+h)) \setind{T_1 > h} | \Upsilon^{(t)}(s) = x, z(s) = \bfz ], \\
	B(h) &= \E[f(\Upsilon^{(t)}(s+h), z(s+h)) \setind{T_1 \leq h < T_2} | \Upsilon^{(t)}(s) = x, z(s) = \bfz ], \\
	C(h) &= \E[f(\Upsilon^{(t)}(s+h), z(s+h)) \setind{T_2 \leq h} | \Upsilon^{(t)}(s) = x, z(s) = \bfz ]. \\
\end{aligned}
\end{equation*}
On the event $\{T_1 > h\}$, it holds that $\Upsilon^{(t)}(s+h) = x$. Further, for $(x, \bfz) \in \overline{\configS}$ and $s \in [0,t]$, let $\Lambda^{(t)}(x, \bfz, s) = \sum_{\bfk, y} \rho^{(t)}_{\bfk,y} (x, \bfz, s)$. It follows that
\begin{equation*}
\begin{aligned}
	\P(T_1 > h | \Upsilon^{(t)}(s) = x, z(s) = \bfz) & = \int_h^{+\infty} \Lambda^{(t)}(x, z(s + t_1), s+t_1) e^{-\int_0^{t_1} \Lambda^{(t)}(x, z(s + u), s + u)du } dt_1 \\
	& = e^{-\int_0^{h} \Lambda^{(t)}(x, z(s + u), s + u)du},
\end{aligned}
\end{equation*}
from which we deduce by the chain rule that 
\begin{equation*}
\begin{aligned}
	\lim_{h \to 0+} \frac{1}{h}&(A(h) - \P(T_1 > h | \Upsilon^{(t)}(s) = x, z(s) = \bfz)f(x, \bfz)) \\
	&  = 
	\lim_{h \to 0+} \P(T_1 > h | \Upsilon^{(t)}(s)  = x, z(s) = \bfz)(f(x, z(s+h)) - f(x, z(s)) \\
	&= \sum_{\bfk, y} \bfz_y \tau_{\bfk}(y,\bfz) \angles{\nabla_\bfz f(x, \bfz)}{\bfk-e(y)}. 
\end{aligned}
\end{equation*}

Let us turn to the case $\{T_1 \leq h < T_2\}$. For $\bfk \in \N^d$ and $y \in \X$, we can compute the probability that the spine is replaced by offspring $\bfk$ and becomes of type $y$ at $T_1$, and does not reproduce anymore until time $s+h$:
\begin{equation*}	
P^{(t)}_{\bfk, y}(h) = \int_0^h \rho^{(t)}(x, z(s+t_1), s+t_1) e^{-\int_0^{t_1} \Lambda^{(t)}(x, z(s+u), s+u)du} e^{-\int_{t_1}^{h} \Lambda^{(t)}(y, z(s+u), s+u)du} dt_1.
\end{equation*}
Hence
\begin{equation*}
\P(T_1 \leq h < T_2 | \Upsilon^{(t)}(s)  = x, z(s) = \bfz) = \sum_{\bfk, y} P^{(t)}_{\bfk, y}(h),
\end{equation*}
and 
\begin{equation*}
B(h) = \sum_{\bfk, y} P^{(t)}_{\bfk, y}(h) f(y, z(s+h)).
\end{equation*}
In particular, $P^{(t)}_{\bfk, y}$ satisfies
\begin{equation*}
	\lim_{h \to 0+} \frac{1}{h} P^{(t)}_{\bfk, y}(h) = \rho^{(t)}_{\bfk, y}(x, \bfz, s).
\end{equation*}
Thus, by continuity of $z$,
\begin{equation*}
	\lim_{h \to 0+} \frac{1}{h}(B(h) - \P(T_1 \leq h < T_2 | \Upsilon^{(t)}(s)  = x, z(s) = \bfz)f(x, \bfz))= \sum_{\bfk, y}  \rho^{(t)}_{\bfk, y}(x, \bfz, s) (f(y, \bfz) - f(x, \bfz)).
\end{equation*}

Finally, a similar computation yields that there exists $C > 0$ such that
\begin{equation*}
\P(T_2 \leq h | \Upsilon^{(t)}(s) = x, z(s) = \bfz) \leq C h^2.
\end{equation*}
By continuity of $z$, given that $z(s) = \bfz$ there exists a compact subset $A \subset \R_+^D$ such that for any $h \in [0,1]$, $(\Upsilon^{(t)}(s+h), z(s+h)) \in \X \times A$. As $f \in \mathcal{C}^1(\overline{\configS})$, it is bounded by some constant $\alpha$ on $\X \times A$, which implies 
\begin{equation*}
\begin{aligned}
	| C(h) -  \P(T_2 \leq h | \Upsilon^{(t)}(s) = x, z(s) = \bfz)f(x, \bfz) | & \leq  2\|f\|_\infty \P(T_2 \leq h | \Upsilon^{(t)}(s) = x, z(s) = \bfz) \\
	& \leq 2 \alpha C h^2 \xrightarrow[h \to 0+]{}0.
\end{aligned}
\end{equation*}
This concludes the proof.

\section{Proof of Lemma \ref{lemma:mups}} 
\label{appdx:mups}

In order to establish Lemma \ref{lemma:mups}, we proceed in two steps. First, we will show that for any $x \in \X$ and $t \geq 0$, $\bfz \mapsto \m(x, \bfz, t)$ is Lipschitz continuous. Rademacher's theorem then ensures that this function is differentiable almost everywhere. Second, we compute $\partial_t \m(x, \bfz, t)$ to show that $\m$ indeed solves Equation \eqref{eq:pde-m}. 

In order to achieve Lipschitz continuity of $\bfz \mapsto \m(x, \bfz, t)$, we want to control the difference between trajectories of $(\Upsilon, z)$ whose population composition starts in different initial conditions. We thus first control the distance between two solutions of dynamical system \eqref{eq:ode}. Throughout the following, we write $\phi(t, \bfz)$ for the solution to Equation \eqref{eq:ode} evaluated at time $t$. 

\begin{lemma}
\label{lemma:semi-flow-lipschitz}
Let $t \geq 0$ and $\mathbf{Z} \subset (0, +\infty)^D$ compact.  Under Assumption \ref{hyp:lipschitz}, for any $\bfz_1, \bfz_2 \in \mathbf{Z}$, there exists $C_\mathbf{Z}(t) > 0$ such that  
\begin{equation*}
 \sup_{s \in [0,t]} \|\phi(s, \bfz_1) - \phi(s, \bfz_2) \|_1 \leq C_\mathbf{Z}(t) \| \bfz_1 - \bfz_2 \|.
\end{equation*}
\end{lemma}

The result is classical for dynamical systems, and is included for completeness.

\begin{proof}
By definition, for any $\bfz \in \mathbf{Z}$,
\begin{equation*}
\phi(t, \bfz) = \bfz + \int_0^t \phi(\bfz, s) A(\phi(\bfz, s)) ds.
\end{equation*}
For any $t \geq 0$ and $\bfz_1, \bfz_2 \in \mathbf{Z}$, there exists a compact subset $\mathbf{K} \subset \R_+^D$ which contains the trajectories $\phi(\bfz_1, \cdot)$ and $\phi(\bfz_2, \cdot)$ on $[0,t]$. 
As $A$ is locally Lipschitz continuous according to Assumption \ref{hyp:lipschitz}, it follows that $\bfz \mapsto \bfz A(\bfz)$ is Lipschitz continuous on $\mathbf{K}$. Thus, there exists $L >0$ such that 
\begin{equation*}
\| \phi(s, \bfz_1) - \phi(s, \bfz_2) \|_1 \leq \|\bfz_1 - \bfz_2 \|_1 + L \int_0^s \| \phi(u, \bfz_1) - \phi(u, \bfz_2) \|_1 du.
\end{equation*}
Hence 
\begin{equation*}
\sup_{s \in [0,t]} \| \phi(s, \bfz_1) - \phi(s, \bfz_2) \|_1 \leq \|\bfz_1 - \bfz_2 \|_1 + L \int_0^t \sup_{\sigma \in [0, u]} \| \phi(\sigma, \bfz_1) - \phi(\sigma, \bfz_2) \|_1 du.
\end{equation*}
The conclusion follows from Gronwall's lemma.
\end{proof}

We may now establish Lipschitz continuity of $\m$. 
\begin{lemma}
\label{lemma:m-lip}
	Under Assumption \ref{hyp:lipschitz}, for any $x \in \X$, the map $(\bfz, t) \mapsto \m(x, \bfz, t)$ is locally Lipschitz continuous on $(0, +\infty)^D \times \R_+$.
\end{lemma}

\begin{proof}
	Let $T > 0$, $x \in \X$ and $\bfz_1, \bfz_2 \in \mathbf{Z} \subset (0, +\infty)^D$ compact. We apply Lemma \ref{lemma:error} to $\zeta_1 = \phi(\cdot, \bfz_1)$ and $\zeta_2 = \phi(\cdot, \bfz_2)$ which according to Lemma \ref{lemma:semi-flow-lipschitz} satisfy Equation \eqref{eq:hyp-zeta} with $\varepsilon_K = C_\mathbf{Z}(t) \|\bfz_1 - \bfz_2 \|$ and $\alpha_K = 0$. In this case, for $j \in \{1,2\}$, $(Y_j, \zeta_j) = (\Upsilon, z)$ with initial condition $\Upsilon(0) = y_0$ and $z(0) = \bfz_j$. For clarity, for $j \in \{1,2\}$, we thus write $(\Upsilon_j, z_j)$ for the spinal constructions with initial condition $(x, \bfz_j)$ defined on the same probability space such that
\begin{equation*}
	\P(\exists s \in [0,t], \Upsilon_1(s) \neq \Upsilon_2(s)) \leq C_\mathbf{Z}(t) \| \bfz_1 - \bfz_2\|_1.
\end{equation*}

Throughout the following, $C_\mathbf{Z}(t)$ is a positive constant whose value may vary from line to line, but which only depends on $t$ and $\mathbf{Z}$. For $t_1, t_2 \in [0,T]$ and $\bfz_1, \bfz_2 \in \mathbf{Z}$, we thus obtain:
\begin{equation*}
	\begin{aligned}
		|\m(x, \bfz_1, t_1) & - \m(x, \bfz_2, t_2)| = \left| \E \left[e^{\int_0^{t_1} \lambda(\Upsilon_1(s), z_1(s))ds} - e^{ \int_0^{t_2} \lambda(\Upsilon_2(s), z_2(s))ds } \right] \right| \\
		& \leq \left| \E \left[e^{\int_0^{t_1} \lambda(\Upsilon_1(s), z_1(s))ds} - e^{\int_0^{t_2} \lambda(\Upsilon_1(s), z_2(s))ds } \right] \right| \\
		& \quad + C_\mathbf{Z}(t) \P(\exists s \in [0,t], \Upsilon_1(s) \neq \Upsilon_2(s)) \\
		& \leq \E \left[ \left| e^{\int_0^{t_1} \lambda(\Upsilon_1(s), z_1(s))ds} - e^{\int_0^{t_2} \lambda(\Upsilon_1(s), z_1(s))ds } \right|  \right]  \\ 
		& + \E \left[ \left| e^{ \int_0^{t_2} \lambda(\Upsilon_1(s), z_1(s))ds } - e^{ \int_0^{t_2} \lambda(\Upsilon_1(s), z_2(s))ds } \right|  \right] \\
		& + C_\mathbf{Z}(T) \|\bfz_1 - \bfz_2 \|_1.
	\end{aligned}
\end{equation*}
where the last inequality follows from Lemma \ref{lemma:semi-flow-lipschitz}. Recall that $\lambda$ is bounded. Thus $t \mapsto \exp \left(\int_0^{t} \lambda(\Upsilon_1(s), z_1(s))ds \right)$ is Lipschitz continuous on $[0,T]$, for a Lipschitz constant which can be chosen to be independent of $x$. Hence 
\begin{equation*}
	\E \left[ \left| \exp \left(\int_0^{t_1} \lambda(\Upsilon_1(s), z_1(s))ds \right) - \exp \left(\int_0^{t_2} \lambda(\Upsilon_1(s), z_1(s))ds \right)\right|  \right] \leq C_\mathbf{Z}(T) | t_2 - t_1|.
\end{equation*}
Further, by Lipschitz continuity of $\lambda$, the following inequality holds almost surely:
\begin{equation*}
	\begin{aligned}
 		\Big| & e^{\int_0^{t_2} \lambda(\Upsilon_1(s), z_1(s))ds } - e^{\int_0^{t_2} \lambda(\Upsilon_1(s), z_2(s))ds}\Big| \\
 		& \leq C_\mathbf{Z}(T) \left| \int_0^{t_2} \lambda(\Upsilon_1(s), z_1(s))ds - \int_0^{t_2} \lambda(\Upsilon_1(s), z_2(s))ds \right| \\
		& \leq C_\mathbf{Z}(T) \sup_{s \in [0,T]} \| z_1(s) - z_2(s)\|_1 \\
 	& \leq C_\mathbf{Z}(T) \| \bfz_1 - \bfz_2\|.
	\end{aligned}
\end{equation*}
This concludes the proof.

\end{proof}

We finally establish the desired result. 
\begin{proof}[Proof of Lemma \ref{lemma:mups}]
	Lemma \ref{lemma:m-lip} ensures by Rademacher's theorem that for any $x \in \X$, the map $\bfz \mapsto \m(x, \bfz, t)$ is differentiable for almost every $t \geq 0$ and $\bfz \in \R_+^D \setminus \{0\}$. We thus consider $t, \bfz$ such that $\m$ is differentiable at $(x, \bfz, t)$, for any $x \in \X$. 
	
Start by noticing that for $h > 0$,
\begin{equation*}
\begin{aligned}
	\m(x, \bfz, t+h) &= \spineE\left[\mathcal{W}(h) \m(\Upsilon(h), z(h), t)\right].
\end{aligned}
\end{equation*}

Let $T_1$ and $T_2$ denote the times of the first and second reproduction event, respectively. Then 
\begin{equation*}
	\begin{aligned}
		\m(x, \bfz, t+h) &= \spineE\left[\mathcal{W}(h)\setind{T_1 > h} \right] \m(x, z(h), t) \\
		&+ \spineE\left[\mathcal{W}(h) \m(\Upsilon(h), z(h), t) \setind{T_1 \leq h < T_2} \right] \\
		&+ \spineE\left[\mathcal{W}(h) \m(\Upsilon(h), z(h), t) \setind{T_2 \leq h} \right].
	\end{aligned}
\end{equation*}
Since $\m$ is differentiable in $\bfz$ at $(t, \bfz)$, we may now proceed as in the Proof of Proposition \ref{prop:lln-spine} to show that 
\begin{equation*}
\lim_{h \to 0+} \frac{1}{h} (\m(x, \bfz, t+h) - \m(x, \bfz, t)) = \mathcal{G}(\m(\cdot, t))(x, \bfz). 
\end{equation*}
This concludes the proof.
\end{proof}

\end{appendix}

\end{document}